\documentclass[12pt,reqno]{amsart}
\usepackage{amsmath,amssymb,amsthm}
\usepackage{mathrsfs}
\usepackage{stmaryrd}

\newtheorem{thm}{Theorem}[section]
\newtheorem{cor}[thm]{Corollary}
\newtheorem{lem}[thm]{Lemma}
\newtheorem{prop}[thm]{Proposition}
\theoremstyle{definition}
\newtheorem{defn}[thm]{Definition}
\newtheorem{nota}[thm]{Notation}
\theoremstyle{remark}

\theoremstyle{remark}

\numberwithin{equation}{section}
\date{December 17, 2008.}
\begin{document}
\title[Moderate deviations for $\mathbb{H}$-valued stationary sequences]{Moderate deviations for stationary sequences of Hilbert-valued bounded random variables}
\author{Sophie Dede}
\address{LPMA, UPMC Universit\'e Paris 6, Case courrier 188, 4, Place Jussieu, 75252 Paris Cedex 05, France.}
\email{sophie.dede@upmc.fr}
\maketitle

\begin{abstract}
In this paper, we derive the Moderate Deviation Principle for stationary sequences of bounded random variables with values in a Hilbert space.
The conditions obtained are expressed in terms of martingale-type conditions. The main tools are martingale approximations and a new Hoeffding inequality for non-adapted sequences of Hilbert-valued random variables. Applications to
Cram\'er-Von Mises statistics, functions of linear processes and stable Markov chains are given.
\end{abstract}
\section{Introduction}
Let $\mathbb{H}$ be a separable Hilbert space with norm $\|.\|_{\mathbb{H}}$
generated by an inner product, $<.,.>_{\mathbb{H}}$, and $(e_l)_{l \geq 1}$ be
an orthonormal basis of $\mathbb{H}$. \\
\indent For the stationary sequence $(X_i)_{i \in \mathbb{Z}}$, of centered random variables with values in $ \mathbb{H}$,
define the partial sums and the normalized process $\{ Z_n(t): t \in [0,1] \}$ by
\begin{equation*}
S_n=\sum_{j=1}^n X_j \ \ \mbox{and} \ \ Z_n(t)=\frac{1}{\sqrt{n}} \sum_{i=1}^{[nt]} X_i + \frac{1}{\sqrt{n}} (nt-[nt]) X_{[nt]+1},
\end{equation*}
$[.]$ denoting the integer part. \\
\indent In this paper, we are concerned with the
Moderate Deviation Principle, for the process $Z_n(.)$, considered as
an element of $C_{\mathbb{H}}([0,1])$, the set of all continuous functions from $[0,1]$ to $\mathbb{H}$. This is a separable
Banach space under the sup-norm $\|x\|_{\infty}=\sup \{ \|x(t)\|_{\mathbb{H}}: t \in[0,1] \} $.
 More generally, we say that a family of random
variables $ \{ Z_n, n >0 \} $ satisfies the Moderate Deviation
Principle (MDP) in $E$, a separable metric space, with speed $a_n \rightarrow 0$,
and good rate function $I(.)$, if the level
sets $\{ x, I(x) \leq \alpha \} $ are compact for all $ \alpha < \infty$, and for all Borel sets $\Gamma$ of $E$,
\begin{align} 
- \inf \{ I(x);x \in \overset{\circ}{\Gamma} \} & \leq & \underset{n \longrightarrow \infty}{\liminf} \ a_n \ \log \ \mathbb{P}(\sqrt{a_n} Z_n \in \Gamma) & & \notag \\
& \leq & \underset{n \longrightarrow \infty}{\limsup} \ a_n \ \log \ \mathbb{P}(\sqrt{a_n} Z_n \in \Gamma) & \leq & - \inf \{ I(x) ; x \in \bar{\Gamma} \}. \label{eqnmdpdef1}
\end{align}
\indent From now, we assume that the stationary sequence $(X_i)_{i
\in \mathbb{Z}}$ is given by $X_i=X_0 \circ T^i$, where $T: \Omega
\longmapsto \Omega$ is a bijective bimeasurable transformation
preserving the probability $\mathbb{P}$ on $(\Omega, \mathcal{A})$.
For a subfield $\mathcal{F}_0$ satisfying $\mathcal{F}_0 \subseteq
T^{-1}(\mathcal{F}_0)$, let $\mathcal{F}_i=T^{-i}(\mathcal{F}_0)$.
By $\|\|X\|_{\mathbb{H}}\|_{\infty}$, we denote the $\mathbf{L}_{\mathbb{H}}^{\infty}$-norm, that
is the smallest $u$ such that $\mathbb{P}(\| X \|_{\mathbb{H}} >u)=0$. \\
\indent When $\mathbb{H}=\mathbb{R}$, Dedecker, Merlev\`ede, Peligrad and
Utev \cite{dedeckermerlevedepeligradutev}
have recently proved (see their Theorem $1$), by using a martingale approximation approach, that:
\begin{thm} \label{theointro}
Assume that $\|X_0\|_{\infty}<\infty$, and that $X_0$ is $\mathcal{F}_0$-measurable. In addition, assume that
\begin{equation*}
\sum_{n=1}^{\infty} \frac{\|\mathbb{E}(S_n \mid \mathcal{F}_0)\|_{\infty}}{n^{3/2}}< \infty,
\end{equation*}
and that there exists $\sigma^2 \geq 0$ with
\begin{equation*}
\underset{n\longrightarrow \infty}{\lim} \Big \|\mathbb{E}\Big (\frac{S_n^2}{n} \Big | \mathcal{F}_0\Big)-\sigma^2 \Big \|_{\infty}=0.
\end{equation*}
Then, for all positive sequences $a_n \rightarrow 0$ and $na_n\rightarrow \infty$, the normalized process $Z_n(.)$ satisfies
the MDP in $C_{\mathbb{R}}([0,1])$, with the good rate function $I_{\sigma}(.)$ defined by
\begin{equation*}
I_{\sigma}(h)=\frac{1}{2 \sigma^2} \int_0^1 \big(h'(u) \big)^2 \, du
\end{equation*}
if simultaneously $\sigma>0$, $h(0)=0$ and $h$ is absolutely continuous, and $I_{\sigma}(h)=\infty$ otherwise.
\end{thm}
The first aim of this paper is to extend the above result to random variables taking their values in a real and separable Hilbert space $\mathbb{H}$.
Indeed, having asymptotic results concerning dependent random variables with values in $\mathbb{H}$ allows for instance, to derive
the corresponding asymptotic results for statistics of the type $\int_0^1 \mid \mathbb{F}_n(t)-\mathbb{F}(t) \mid^2 \, \mu(dt)$ where $\mathbb{F}(.)$ is the cumulative distribution function
of a strictly stationary sequence of real random variables $(Y_i)_{i \in \mathbb{Z}}$ and $\mathbb{F}_n(.)$ is the corresponding empirical distribution function
(see Section 3.4). \\
\indent On an other hand, since Theorem \ref{theointro} is stated for adapted sequences, the second aim of this paper is to extend this result to
non-adapted sequences. \\
\indent To extend Theorem \ref{theointro} to non-adapted sequences of Hilbert-valued random variables, we use a similar
martingale approach as done for instance in Voln\'y \cite{volny} for the central limit theorem.
In infinite dimensional cases, the authors have essentially considered i.i.d or triangular arrays of i.i.d random variables
(see for instance de Acosta \cite{deacosta}, Borovkov and Mogulskii \cite{borovkovmogulskii1} \cite{borovkovmogulskii2}, Ledoux \cite{ledoux}, ...).
However for dependent sequences with values in functional spaces, there are few results available in the literature.
Since our approach is based on martingale approximation, we first extend Puhalskii \cite{puhalskii} results for
$\mathbb{R}^d$-valued martingale differences sequences to the $\mathbb{H}$-valued case (see Section 4.2).
In Section 2.1, we derive a Hoeffding inequality for a sequence of
non-adapted Hilbert-valued random variables.
Section 4 is dedicated to the proofs.

\section{Main Results}
We begin with some notations,
\begin{nota} \label{Notation2}
 For any real $p \geq 1$, denote by
$\mathbf{L}^p_{\mathbb{H}}$ the space of $\mathbb{H}$-valued random
variables $X$ such that
$\|X\|^p_{\mathbf{L}^p_{\mathbb{H}}}=\mathbb{E}(\|X\|^p_{\mathbb{H}})$
is finite. For example, $\mathbf{L}_{\mathbb{H}}^1([0,1])$ is the space of $\mathbb{H}$-valued Bochner integrable functions on $[0,1]$.
\end{nota}
\subsection{A Hoeffding inequality}
$ \newline $
\indent Firstly, we start by establishing a maximal inequality, which is obtained through a generalization of the ideas in Peligrad,
Utev and Wu \cite{peligradutevwu}.
\begin{thm} \label{theohoeffding}
Assume that $\|\|X_0\|_{\mathbb{H}}\|_{\infty} < \infty$. For any $x>0$, we have
\begin{equation} \label{eqnhoeffding1}
\mathbb{P} \Big( \underset{1 \leq i \leq n}{\max} \|S_i\|_{\mathbb{H}} \geq x \Big) \leq 2 \sqrt{e} \exp \left(- \frac{x^2}{4 n \big( \|\|X_0\|_{\mathbb{H}}\|_{\infty} + C \Delta \big)^2} \right),
\end{equation}
for some constant $C >0$ and
\begin{equation*} \label{eqnhoeffding2}
\Delta=\sum_{j=1}^n \frac{1}{j^{3/2}}\big( \|\|\mathbb{E}(S_j \mid \mathcal{F}_0) \|_{\mathbb{H}} \|_{\infty} + \|\| S_j - \mathbb{E}(S_j \mid \mathcal{F}_j)\|_{\mathbb{H}} \|_{\infty} \big).
\end{equation*}
\end{thm}

\subsection{The Moderate Deviation Principle}
$\newline$
 \indent Before establishing our main result, we need more definitions.
\begin{defn} \label{deftrace}
A nonnegative self-adjoint operator $ \Gamma $ on $\mathbb{H}$ will be called an $\mathcal{S}(\mathbb{H})$-operator, if it
has finite trace, i.e, for some ( and therefore every) orthonormal basis $(e_l)_{l \geq 1}$ of $\mathbb{H}$, $\sum_{l \geq 1} < \Gamma e_l,e_l>_{\mathbb{H}}< \infty$.
\end{defn}
Let
\begin{align*}
 \mathcal{AC}_0([0,1])& =\{ \phi \in C_{\mathbb{H}}([0,1]):  \ \mbox{there exists} \ g \in \mathbf{L}_{\mathbb{H}}^1([0,1]) \\
                      & \qquad \qquad \ \mbox{such that} \ \phi(t)=\int_0^t g(s) \, ds \ \mbox{for} \ t \in [0,1] \}.
\end{align*}
Now, we give the extension of Theorem \ref{theointro}.
\begin{thm} \label{theomdpfonctionnel}
Assume that $\|\|X_0 \|_{\mathbb{H}} \|_{\infty} < \infty $. Moreover, assume that
\begin{equation} \label{conditionfonctionnel}
\sum_{n \geq 1 } \frac{1}{n^{3/2}} \|\|\mathbb{E}(S_n \mid \mathcal{F}_0) \|_{\mathbb{H}} \|_{\infty} < \infty \ \ \mbox{and} \ \ \sum_{n \geq 1} \frac{1}{n^{3/2}} \|\|S_n- \mathbb{E}(S_n \mid \mathcal{F}_n) \|_{\mathbb{H}}\|_{\infty} < \infty,
\end{equation}
and that there exists $Q \in \mathcal{S}(\mathbb{H})$ such that
\begin{enumerate}
\item[$i.$] for all $k$, $l$ in $\mathbb{N}^{\ast}$,
\begin{equation} \label{eqnfonct1}
\underset{n \longrightarrow \infty}{\lim} \Big \| \frac{1}{n} \mathbb{E}(<S_n,e_k>_{\mathbb{H}} <S_n,e_l>_{\mathbb{H}} \mid \mathcal{F}_0)- <Q e_k,e_l>_{\mathbb{H}} \Big \|_{\infty}=0,
\end{equation}
\item[$ii.$]
\begin{equation} \label{eqnfonct2}
\underset{ n \longrightarrow \infty} {\lim} \Big \| \frac{1}{n} \mathbb{E}(\|S_n \|_{\mathbb{H}}^2 \mid \mathcal{F}_0) - \mbox{Tr}(Q) \Big \|_{\infty}=0.
\end{equation}
\end{enumerate}
Then, for all positive sequences $a_n$ with $a_n \rightarrow 0$ and $n a_n \rightarrow \infty$, the
process $ Z_n(.)$ satisfies the functional MDP in $C_{\mathbb{H}}([0,1])$ with the good
rate function,
\begin{equation} \label{eqnfonctiontaux}
 I(\phi)=\left \lbrace \begin{array} {l}
                    \int_0^1 \Lambda^{\ast}(\phi'(t)) \, dt \ \ \ \mbox{if} \ \phi \in \mathcal{AC}_0([0,1])  \\ \\
                    + \infty \ \ \mbox{otherwise,}
                    \end{array}
                 \right.
 \end{equation}
where $\Lambda^{\ast}$ is given by:
\begin{equation} \label{eqnlambdaetoile}
 \Lambda^{\ast}(x)=\underset{y \in \mathbb{H}}{\sup} \big (<y,x>_{\mathbb{H}} -\frac{1}{2} <y,Qy>_{\mathbb{H}} \big).
\end{equation}
\end{thm}
As an immediate consequence, we have:
\begin{cor} \label{theomdp}
Under the same notations and assumptions of  Theorem \ref{theomdpfonctionnel}, we have that, for all positive sequences $a_n$
with $a_n \rightarrow  0$ and $n a_n \rightarrow \infty$,
$n^{-1/2} S_n$ satisfies the MDP in $\mathbb{H}$  with the good rate
function, $\Lambda^{\ast}$ defined in (\ref{eqnlambdaetoile}).
\end{cor}
Since $\mbox{Tr}(Q)< \infty$,  $Q$ is a compact operator. If $x \in Q(\mathbb{H})$, then there is $z \in \mathbb{H}$, such that $x=Qz$. Hence, the rate function is
$\nobreak$
\begin{equation*}
\forall \ x \in Q(\mathbb{H}), \ \Lambda^{\ast}(x)=\frac{1}{2} <z,Qz>_{\mathbb{H}} = \frac{1}{2} <z,x>_{\mathbb{H}}.
\end{equation*}
If $x \not \in Q(\mathbb{H})$, we have $\Lambda^{\ast}(x)=+\infty$.
In particular, if $Q$ is injective, $(\lambda_i)_{i \geq 1} $
are its eigenvalues, and $(f_i)_{i \geq 0}$ the associated eigenvectors, we can
simplify the rate function,
\begin{equation*} \label{eqntheomdp2}
\forall  \ x \in Q(\mathbb{H}), \ \Lambda^{\ast}(x)= \frac{1}{2} \sum_{i \geq 1}
\frac{1}{\lambda_i} <x,f_i>^2_{\mathbb{H}}.
\end{equation*}
\indent The following corollary gives simplified conditions for the MDP.
\begin{cor} \label{corollaire}
Assume that $ \| \| X_0 \|_{\mathbb{H}} \|_{\infty} < \infty$. Moreover, assume that
\begin{equation} \label{conditionfonctionnelcor}
\sum_{n \geq 1 } \frac{1}{\sqrt{n}} \|\|\mathbb{E}(X_n \mid \mathcal{F}_0) \|_{\mathbb{H}} \|_{\infty} < \infty \ \ \mbox{and} \ \ \sum_{n \geq 1} \frac{1}{\sqrt{n}} \|\|X_{-n}- \mathbb{E}(X_{-n} \mid \mathcal{F}_0) \|_{\mathbb{H}}\|_{\infty} < \infty,
\end{equation}
and that for all $i, \ j$ in $\mathbb{N}^{\ast}$,
\begin{enumerate}
\item[$1.$] for all  $k$, $l$ in $\mathbb{N}^{\ast}$,
\begin{equation} \label{eqnfonct1cor}
\underset{n \longrightarrow \infty}{\lim} \| \mathbb{E}(<X_i,e_k>_{\mathbb{H}}<X_j,e_l>_{\mathbb{H}} \mid \mathcal{F}_{-n})-  \mathbb{E}(<X_i,e_k>_{\mathbb{H}}<X_j,e_l>_{\mathbb{H}}) \|_{\infty}=0,
\end{equation}
\item[$2.$]
\begin{equation} \label{eqnfonct2cor}
\underset{ n \longrightarrow \infty} {\lim} \|  \mathbb{E}(<X_i,X_j >_{\mathbb{H}} \mid \mathcal{F}_{-n}) - \mathbb{E}(<X_i,X_j >_{\mathbb{H}}) \|_{\infty}=0.
\end{equation}
\end{enumerate}
Then, the conclusion of Theorem \ref{theomdpfonctionnel} holds, with $Q$ defined by
\begin{equation*}
\mbox{for all $k$, $l \geq 1$,} \ \ <e_k,Qe_l>_{\mathbb{H}}=\sum_{p \in \mathbb{Z}} \mathbb{E}(<X_0,e_k>_{\mathbb{H}}<X_p,e_l>_{\mathbb{H}}).
\end{equation*}
\end{cor}
\subsection{Functional law of the iterated logarithm}
$\newline$
Throughout this section, let $ \beta(n)=\sqrt{2n \log \log n}, \ n \geq 3$. Let $\tilde{S_n}(.)$ be the process $\{ \overset{\sim}{S_n}(t)=\sum_{i=1}^{[nt]} X_i +(nt-[nt]) X_{[nt] +1}: t \in [0,1] \}$.
\begin{thm} \label{functionallaw}
Assume that $\|\|X_0\|_{\mathbb{H}}\|_{\infty}< \infty$. Assume in addition that (\ref{conditionfonctionnel}),(\ref{eqnfonct1}) and (\ref{eqnfonct2}) hold.
Then, with probability $1$, the following sequence
\begin{equation*}
\Big \{ \xi_n(.)=\frac{\overset{\sim}{S_n}(.)}{\beta(n)} \Big \}_{n \geq 1}
\end{equation*}
is relatively compact in $C_{\mathbb{H}}([0,1])$ and the set of its limit points is precisely the compact set
\begin{equation*}
\mathcal{K}=\{ \phi \in C_{\mathbb{H}}([0,1]), \ \mbox{such that} \ 2 I(\phi)\leq 1 \}.
\end{equation*}
\end{thm}
\noindent \textit{Proof of Theorem \ref{functionallaw}}. It can be proved by the arguments of the proof of Theorem $3.1$ in Hu and Lee \cite{hulee} ( see also Deuschel and Stroock \cite{deuschelstroock}). $\square$
\section{applications}
\subsection{$\phi$-mixing sequences}
$\newline$
Recall that if $Y$ is a random variable with values in a Polish space $\mathcal{Y}$ and if $\mathcal{F}$ is a $\sigma$-field,
the $\phi$-mixing coefficient between $\mathcal{F}$ and $\sigma(Y)$ is defined by
\begin{equation*}
\phi(\mathcal{F}, \sigma(Y))=\underset{ A \in \mathcal{B}(\mathcal{Y})}{\sup } \| \mathbb{P}_{Y \mid \mathcal{F}} (A)-\mathbb{P}_{Y}(A ) \|_{\infty}.
\end{equation*}
For the sequence $(X_i)_{i \in \mathbb{Z}}$, let
\begin{equation*}
\phi_1 (n)=\phi(\mathcal{F}_0, \sigma(X_n)) \ \mbox{and} \  \phi_{2}(n)=\underset{i > j \geq n}{\sup} \phi (\mathcal{F}_{0}, \sigma(X_i,X_j)).
\end{equation*}
\begin{prop} \label{propmelange}
Assume that $ \|\|X_0\|_{\mathbb{H}}\|_{\infty} < \infty$ and $X_0$ is $\mathcal{F}_0$-measurable. Then, for all $x \geq 0$, we have
\begin{equation} \label{eqnmelange1}
\mathbb{P} \Big( \underset{1 \leq i \leq n}{\max} \|S_i\|_{\mathbb{H}} \geq x \Big) \leq 2 \sqrt{e} \exp \left(- \frac{x^2}{4 n  \|\|X_0\|_{\mathbb{H}}\|^2_{\infty}\big( 1 +6 C \sum_{j \geq 1} j^{-1/2} \phi_1 (j) \big )^2} \right),
\end{equation}
for the same positive constant $C$ defined in Theorem \ref{theohoeffding}.
\end{prop}
\noindent \textit{Proof of Proposition \ref{propmelange}}.
Applying triangle inequality and changing the order of summation, observe that
\begin{equation*}  \label{eqnmelange2}
 \sum_{n \geq 1} \frac{1}{n^{3/2}} \| \| \mathbb{E}(S_n \mid \mathcal{F}_0) \|_{\mathbb{H}} \|_{\infty} \leq 3 \sum_{n \geq 1 } \frac{1}{\sqrt{n}} \|\|\mathbb{E}(X_n \mid \mathcal{F}_0) \|_{\mathbb{H}} \|_{\infty}.
\end{equation*}
Since $\mathbb{E}(X_0)=0$, we have
\begin{equation*}
\|\|\mathbb{E}(X_n \mid \mathcal{F}_0)\|_{\mathbb{H}}\|_{\infty} \leq 2\|\|X_0\|_{\mathbb{H}}\|_{\infty} \phi_1(n).
\end{equation*}
\hfill $\square$ \\
\noindent Next, we have a Moderate Deviation Principle.
\begin{prop} \label{propMDPmelange}
Assume that $\|\|X_0\|_{\mathbb{H}}\|_{\infty} < \infty$ and $X_0$ is $\mathcal{F}_0$-measurable. If
\begin{equation*}
\sum_{n \geq 1} \frac{1}{\sqrt{n}} \phi_1(n) < \infty \ \mbox{and} \ \phi_2(n) \underset{n \longrightarrow \infty}{\longrightarrow} 0,
\end{equation*}
then the conclusion of Theorem \ref{theomdpfonctionnel} holds.
\end{prop}
\subsection{Functions of linear processes}
$\newline$
\indent In this section, we shall focus on functions of $\mathbb{H}$-valued linear processes,
\begin{equation} \label{eqn423}
X_k=f \Big(\sum_{i \in \mathbb{Z}} c_i(\varepsilon_{k-i}) \Big )-\mathbb{E}\Big(f\Big(\sum_{i \in \mathbb{Z}} c_i( \varepsilon_{k-i})\Big)\Big),
\end{equation}
where $f: \mathbb{H} \rightarrow \mathbb{H}$, $(c_i)_{i \in \mathbb{Z}}$ are linear operators from $\mathbb{H}$ to $\mathbb{H}$ and $(\varepsilon_i)_{i \in \mathbb{Z}}$ is
a sequence of i.i.d $\mathbb{H}$-valued random variables such that $\|\|\varepsilon_0\|_{\mathbb{H}}\|_{\infty}< \infty$.\\
\indent The sequence $\{X_k\}_{k \geq 1}$ defined by (\ref{eqn423}) is a natural extension of the multivariate linear processes. These types of processes with values in functional spaces also facilitate the study of
estimating and forecasting problems for several classes of continuous time processes (see Bosq \cite{bosq}).\\
\indent We denote by $\|.\|_{L(\mathbb{H})}$, the operator norm.
We shall give sufficient conditions for the Moderate Deviation Principle in terms of the regularity of the function $f$. \\
\indent Let $\delta(\varepsilon_0)=2 \inf \{ \| \| \varepsilon_0 -x \|_{\mathbb{H}}\|_{\infty}, x \in \mathbb{H} \} $ and define the modulus of continuity of $f$ by
\begin{equation*} \label{eqn424}
w_f(h)=\underset{ \|t\|_{\mathbb{H}} \leq h , \ x \in  \mathbb{H}}{\sup} \|f(x+t)-f(x) \|_{\mathbb{H}}.
\end{equation*}
\begin{prop} \label{propfunclinearprocess}
Assume that
\begin{equation*} \label{eqn426}
\sum_{i \in \mathbb{Z}} \|c_i\|_{L(\mathbb{H})} < \infty,
\end{equation*}
and that $X_k$ is defined as in (\ref{eqn423}). If moreover
\begin{equation} \label{eqn427}
\sum_{n \geq 1} \frac{1}{\sqrt{n}} w_f \Big(\delta(\varepsilon_0) \sum_{k \geq n} \|c_k\|_{L(\mathbb{H})}\Big)+\sum_{n \geq 1} \frac{1}{\sqrt{n}} w_f \Big(\delta(\varepsilon_0) \sum_{k \leq- n} \|c_k\|_{L(\mathbb{H})} \Big)< \infty,
\end{equation}
then the conclusion of Theorem \ref{theomdpfonctionnel} holds.
\end{prop}
In particular, if $\|c_i\|_{L(\mathbb{H})}=O(\rho ^{\mid i \mid}), \ 0 < \rho < 1$, the condition (\ref{eqn427}) is equivalent to
\begin{equation} \label{eqn3001}
\int_0^1 \frac{w_f(t)}{t \sqrt{|\log t|}} \, dt < \infty.
\end{equation}
For example, if $w_f(t) \leq D  |\log t |^{-\gamma}$ for some $D>0$ and some $\gamma>1/2$, then (\ref{eqn3001}) holds.\\ \\
\underline{\textbf{Remark}:} Under a Cramer type condition, Mas and Menneteau
\cite{masmenneteau} were interested in the MDP for the asymptotic
behavior of the empirical mean of $\overline{X_n}=\frac{1}{n}\sum_{k=1}^n
X_k$, where for all $n\geq 1$, $X_n$ is an autoregressive process:
$X_n=\sum_{j=0}^{\infty}\rho^j(\varepsilon_{n-j})$. Here,
$(\varepsilon_k)_{k\in \mathbb{Z}}$ is a sequence of i.i.d
Hilbert-valued centered random variables satisfying a Cramer type condition and $\rho$ is a bounded
Hilbert linear operator, satisfying
$\sum_{j=0}^{\infty}\|\rho^j\|_{\mathbf{L}(\mathbb{H})}<\infty$.
They gave also the MDP for the difference
between the empirical and theoretical covariance operators.

\subsection{Stable Markov chains}
$\newline$
\indent Let $(Y_n)_{n \geq 0}$ be a stationary Markov chain of $\mathbb{H}$-valued bounded random variables.
Denote by
$\mu$ the law of $Y_0$ and by $K$ its transition kernel. \\
 Let
\begin{equation}\label{eqnstablemarkovchain00}
X_k=f(Y_k)-\mathbb{E}(f(Y_k)).
\end{equation}
For all Lipschitz functions $g:\mathbb{H} \longrightarrow \mathbb{H}$, let
\begin{equation*}
\mbox{Lip}(g)=\underset{x,y \in \mathbb{H}}{\sup} \frac{\|g(x)-g(y)\|_{\mathbb{H}}}{\|x-y\|_{\mathbb{H}}}.
\end{equation*}
\indent We write $K(g)$ and $K^n(g)$ respectively for the functions $y \mapsto \int g(z)K(y,dz)$ and $ y \mapsto \int g(z)K^n(y,dz)=\mathbb{E}(g(Y_n) \mid Y_0=y)$.
\begin{prop} \label{propstablemarkovchain} Assume that the transition kernel $K$ satisfies $ \mathrm{Lip}(K^n(g)) \leq C \rho^n \mathrm{Lip}(g)$ for some $\rho<1$ and any Lipschitz function $g$.
If $f$ is Lipschitz and $X_i$ is defined by (\ref{eqnstablemarkovchain00}), then the normalized process $Z_n(.)$ satisfies the MDP in $C_{\mathbb{H}}([0,1])$.
\end{prop}
Let $(Y_n)_{n \geq 0}$ satisfying the equation
$Y_n=F(Y_{n-1}, \xi_n)$ for some measurable map $F$ and some sequence $(\xi_i)_{i \geq 0}$ of i.i.d random variables independent of $Y_0$.
\begin{cor} \label{corstablemarkovchain} Assume that for all $x, \ y \in \mathbb{H}$,  \ $  \|\|F(x,\xi_1)-F(y, \xi_1)\|_{\mathbb{H}}\|_{\mathbf{L}^1(\mathbb{R},\mathbb{P}_{\xi_1})} \leq \rho \|x-y\|_{\mathbb{H}}$ with $\rho<1$. If $f$ is Lipschitz and
$X_i$ is defined by (\ref{eqnstablemarkovchain00}),
then the normalized process $Z_n (.)$ satisfies the MDP in $C_{\mathbb{H}}([0,1])$.
\end{cor}
\noindent \textit{Proof of Corollary \ref{corstablemarkovchain}}.
The condition: for all $x, \ y \in \mathbb{H}$, $\|\|F(x,\xi_1)-F(y, \xi_1)\|_{\mathbb{H}}\|_{\mathbf{L}^1(\mathbb{R},\mathbb{P}_{\xi_1})} \leq \rho \|x-y\|_{\mathbb{H}}$ with $\rho<1$,
implies that, for all $f$ Lipschitz function,
\begin{equation*}
\mbox{Lip}(K(f)) \leq \rho \mbox{Lip}(f).
\end{equation*}
Indeed, for all $y $ in $\mathbb{H}$, we get
\begin{equation*}
K(f)(y)=\mathbb{E}(f(Y_1) \mid Y_0=y)=\mathbb{E}(f(F(y,\xi_1)) \mid Y_0=y)=\int f(F(y,z)) \mathbb{P}_{\xi_1}(dz).
\end{equation*}
Hence, for all $x, \ y \in \mathbb{H}$, we derive
\begin{eqnarray*}
\mid K(f)(x)-K(f)(y) \mid & \leq & \int \mid f(F(x,z))-f(F(y,z)) \mid \mathbb{P}_{\xi_1}(dz) \\
                          & \leq & \mbox{Lip}(f) \|\|F(x,\xi_1)-F(y,\xi_1)\|_{\mathbb{H}}\|_{\mathbf{L}^1(\mathbb{R}, \mathbb{P}_{\xi_1})} \\
                          & \leq & \rho \mbox{Lip}(f) \|x-y\|_{\mathbb{H}}.
\end{eqnarray*}
It follows that
\begin{equation*}
\mbox{Lip}(K(f)) \leq \rho \mbox{Lip}(f),
\end{equation*}
so that
\begin{equation*}
\mbox{Lip}(K^n(f)) \leq \rho^n \mbox{Lip}(f).
\end{equation*}
We conclude by applying Proposition \ref{propstablemarkovchain}.
\hfill $\square$
\subsection{Moderate Deviation Principle for the empirical distribution function in $\mathbf{L}^2$}
$\newline$
\indent Let $Y=(Y_i)_{i \in \mathbb{Z}}$ be a strictly stationary sequence of real-valued random variables with common distribution function $\mathbb{F}$.
Set $\mathcal{F}_0= \sigma(Y_i,i \leq 0)$. We denote $\mathbb{F}_n$, the empirical distribution function of $Y$:
\begin{equation*}
\forall \ t \in \mathbb{R}, \ \mathbb{F}_n(t) = \frac{1}{n} \sum_{i=1}^{n} \mathbf{1}_{Y_i \leq t}.
\end{equation*}
Note that for any probability measure $\mu$ on $\mathbb{R}$, the
random variable $X_i=\{ t \mapsto \mathbf{1}_{Y_i \leq t}
-\mathbb{F}(t): \ t \in \mathbb{R} \}$ may be viewed as a random variable with values in
the Hilbert space $\mathbb{H}:=\mathbf{L}^2(\mathbb{R},\mu)$.
Hence to derive the MDP for $n(\mathbb{F}_n-\mathbb{F})$ we shall apply Corollary \ref{corollaire} to the random variables $(X_i)_{i \geq 1}$.
With this aim, we first recall
the following dependance coefficients
from Dedecker and Prieur \cite{dedeckerprieur07}.
\begin{defn} \label{coeffmelange}
Let $(\Omega, \mathcal{A}, \mathbb{P})$ be a probability space, let $\mathcal{F}$ be a sub $\sigma$-algebra of $\mathcal{A}$. Let
$Y=(Y_1,...,Y_k)$ be a random variable with values in $\mathbb{R}^k$. Let $\mathbb{P}_{Y}$ be the distribution of $Y$ and let
$\mathbb{P}_{Y \mid \mathcal{F}}$ be a conditional distribution of $Y$ given $\mathcal{F}$. For $1 \leq i \leq k$ and
$t$ in $\mathbb{R}$, let $g_{t,i}(x)=\mathbf{1}_{x \leq t} - \mathbb{P}(Y_i \leq t)$. Define the random variable
\begin{equation*}
b(\mathcal{F},Y_1,...,Y_k)=\underset{(t_1,...,t_k) \in \mathbb{R}^k}{\sup} \Bigg | \int \prod_{i=1}^k g_{t_i,i}(x_i) \mathbb{P}_{Y \mid \mathcal{F}}(dx) - \int \prod_{i=1}^k g_{t_i,i}(x_i) \mathbb{P}_Y(dx) \Bigg |
\end{equation*}
with $ \mathbb{P}_{Y \mid \mathcal{F}}(dx)=\mathbb{P}_{Y \mid \mathcal{F}}(dx_1,...,dx_k)$ and $\mathbb{P}_Y(dx)=\mathbb{P}_Y(dx_1,...,dx_k)$. \\
\indent For the stationary sequence $(Y_i)_{i \in \mathbb{Z}}$, define the coefficient $\tilde{\phi_k}$, for any integer $k \geq 1$, by
\begin{equation*}
\tilde{\phi_k}(n)=\underset{1 \leq l \leq k}{\max} \ \underset{i_l>...>i_1 \geq n}{\sup} \|b(\mathcal{F}_0,Y_{i_1},...,Y_{i_l})\|_{\infty}.
\end{equation*}
\end{defn}
\begin{prop} \label{propcramervonmises}
If
\begin{equation} \label{eqncvm}
\sum_{n \geq 1} \frac{1}{\sqrt{n}} \tilde{\phi_2}(n)< \infty,
\end{equation}
then $ \{ \sqrt{n}(\mathbb{F}_n(t)-\mathbb{F}(t)), \ t \in \mathbb{R} \} $
satisfies the MDP in $\mathbf{L}^2(\mathbb{R},\mu)$ with the good rate function,
\begin{equation} \label{eqncvm4}
\forall \ f \in \mathbf{L}^2(\mathbb{R},\mu), \ I(f)=\underset{g \in \mathbf{L}^2(\mathbb{R},\mu)} {\sup} \big (<f,g>_{\mathbf{L}^2(\mathbb{R},\mu)}-\frac{1}{2}<g,Qg>_{\mathbf{L}^2(\mathbb{R},\mu)} \big ),
\end{equation}
where $Q$ is defined as follows,
for all $(f,g)$ in $\mathbf{L}^2(\mathbb{R},\mu) \times \mathbf{L}^2(\mathbb{R},\mu)$,
\begin{equation*}
Q(f,g) = \int_{\mathbb{R}^2} f(s)g(t) C(s,t) \, \mu(dt) \mu(ds)
\end{equation*}
with
\begin{equation*}
C(s,t) = \mathbb{F}(t \wedge s)-\mathbb{F}(t) \mathbb{F}(s) + 2 \sum_{k \geq 1} ( \mathbb{P}(Y_0 \leq t, Y_k \leq s )-\mathbb{F}(t) \mathbb{F}(s)).
\end{equation*}
\end{prop}
\indent If we use the contraction principle in Dembo and Zeitouni \cite{dembozeitouni}, with the continuous function $f:x \longmapsto \|x\|_{\mathbf{L}^2(\mathbb{R},\mu)}$, the
 Cram\'er-Von Mises statistics,
 \begin{equation*}
 \sqrt{n} \big (\int_{\mathbb{R}}(\mathbb{F}_n(t)-\mathbb{F}(t))^2 \, \mu(dt) \big)^{1/2}
 \end{equation*}
 satisfies the MDP in $\mathbb{R}$ with the good rate function,
 \begin{equation*}
 \forall \ y \geq 0, \ I'(y)=\frac{1}{2} \frac{1}{\nu} y^2,
 \end{equation*}
 where $\nu=\underset{k}{\max} (\lambda_k)$, the $\lambda_k$'s are the eigenvalues of the covariance function $Q$. \\
\indent Now we suppose that $Y_i \in [0,1]$ and $ \mu(dt)=dt$. Always, by the contraction principle in Dembo and Zeitouni \cite{dembozeitouni} with the continuous function,
\begin{eqnarray*}
u:  \mathbf{L}^2([0,1],\mu) & \longrightarrow & \mathbb{R}^{+} \\
 g & \longmapsto & \|g\|_{\mathbf{L}^1([0,1],\mu)},
\end{eqnarray*}
we prove that the Kantorovitch distance $  \sqrt{n} \|\mathbb{F}_n-\mathbb{F}\|_{\mathbf{L}^1([0,1],\mu)}$
satisfies the MDP in $ \mathbb{R}^{+}$, with the good rate function,
\begin{equation*}
\forall \ y \in \mathbb{R}^{+}, \ J(y)=\inf \{ I(f), f \in \mathbf{L}^2([0,1],\mu), y = \|f\|_{\mathbf{L}^1([0,1],\mu)} \}.
\end{equation*}
We deduce from the proof of Theorem in Ledoux \cite{ledoux} page $274$, that
\begin{equation*}
\forall  \ y \geq 0, \ J(y)=\frac{1}{2} \frac{y^2}{\sigma(Z)^2},
\end{equation*}
where $Z$ is a random variable with covariance function $Q$ defined in Proposition \ref{propcramervonmises} and $\sigma(Z)=\underset{ \|g\|_{\infty} \leq 1}{\sup} \Big(\mathbb{E}\Big(\int_0^1 g(t)Z(t) \, dt\Big)^2\Big)^{1/2} $. \\
By using a remark $(8.22)$ in Ledoux and Talagrand \cite{ledouxtalagrand} page $216$, we also have
\begin{equation*}
\underset{n \longrightarrow \infty}{\limsup} \frac{\sqrt{n}\|\mathbb{F}_n-\mathbb{F}\|_{\mathbf{L}^1([0,1],\mu)}}{\sqrt{2n\log \log n}}=\sigma(Z).
\end{equation*}
\underline{\textbf{Remark}:} $(\ref{eqncvm})$ is satisfied for a large class of dependent sequences. For instance, it is verified for the class of
expanding maps as considered in Dedecker and Prieur \cite{dedeckerprieur07}.
\section{Proofs}
\subsection{Hoeffding inequality's proof}
\subsubsection{Technical propositions}
$\newline$
 The proofs of the following lemma and propositions use the same ideas as in the proof of Theorem $1$ in Peligrad, Utev and Wu \cite{peligradutevwu}
(see also Mackey and Tyran-Kami\'nska \cite{mackeytyran}).
\begin{lem} \label{lemhoeffdingdiffmart} Let $\{ Z_k \}_{k \in \mathbb{Z}}$ be a stationary sequence of martingale differences with values in $\mathbb{H}$. For all integer $n \geq 1$ and  real $p \geq 1$,
we have
\begin{equation} \label{eqnlemhoeffdingdiffmart1}
\mathbb{E} \Big(\underset{1 \leq i \leq n}{\max} \|Z_1+...+Z_i\|_{\mathbb{H}}^{2p} \Big) \leq 2^{p+1} \Gamma(p+1) \  n^p \|\|Z_1\|_{\mathbb{H}}\|_{\infty}^{2p}.
\end{equation}
\end{lem}
\noindent \textit{Proof of Lemma \ref{lemhoeffdingdiffmart}}.
We note here, $S_i=\sum_{k=1}^i Z_k$. Applying an inequality given in Pinelis \cite{pinelis94} (Theorem 3.5) and using stationarity, we have
\begin{eqnarray*}
\mathbb{E}\Big(\underset{1 \leq i \leq n}{\max} \|S_i\|_{\mathbb{H}}^{2p} \Big) & = & \int_0^{\infty} (2p) z^{2p-1} \mathbb{P}\Big(\underset{1 \leq i \leq n}{\max} \|S_i\|_{\mathbb{H}} \geq z \Big) \, dz \notag \\
                                                            & \leq & 2 \int_0^{\infty} (2p)  z^{2p-1} \exp \Big(-\frac{z^2}{2n \|\|Z_1\|_{\mathbb{H}}\|_{\infty}^2} \Big) \, dz. \label{eqnlemhoeffdingdiffmart2}
\end{eqnarray*}
By the change of variable $u=\frac{z}{\sqrt{n} \|\|Z_1\|_{\mathbb{H}}\|_{\infty}}$,
\begin{equation*}\label{eqnlemhoeffdingdiffmart3}
\int_0^{\infty} z^{2p-1} \exp \Big(-\frac{z^2}{2(\sqrt{ n} \|\|Z_1\|_{\mathbb{H}} \|_{\infty})^2}\Big) \, dz =  n^p \|\|Z_1\|_{\mathbb{H}}\|_{\infty}^{2p}
\int_0^{\infty} u^{2p-1} \exp \big(-\frac{u^2}{2} \big) \, du.
\end{equation*}
Next
\begin{equation*}
\int_0^{\infty} u^{2p-1} \exp \big (-\frac{u^2}{2} \big ) \, du = 2^{p-1} \Gamma(p).
\end{equation*}
Therefore, we conclude that
\begin{equation*} \label{eqnlemhoeffdingdiffmart4}
\mathbb{E} \Big (\underset{1 \leq i \leq n} {\max} \|S_i\|_{\mathbb{H}}^{2p} \Big ) \leq 2^{p+1} \Gamma(p+1) \ n^p \|\|Z_1 \|_{\mathbb{H}} \|_{\infty}^{2p}.
\end{equation*}
\hfill $\square$\\
\indent The next proposition is a generalization of Lemma \ref{lemhoeffdingdiffmart} to an adapted stationary sequence.
\begin{prop} \label{lemhoeffdingstatprocadap}
Let $n$, $q$ be integers such that $n \geq 1$, $2^{q-1} \leq n < 2^q$. Assume that $ \|\|Z_0 \|_{\mathbb{H}} \|_{\infty} < \infty$,
and $Z_0$ is $\mathcal{F}_0$-measurable. Let $Z_i=(Z_0 \circ T^i)_{i \in \mathbb{Z}}$. Then, for all real $p \geq 1$,
\begin{equation} \label{eqnlemhoeffdingstatprocadap1}
\mathbb{E} \Big(\underset{1 \leq i \leq n}{\max} \Big \| \sum_{j=1}^i Z_j \Big \|_{\mathbb{H}}^{2p} \Big) ^{1/2p} \leq (2^{p+1} \Gamma(p+1))^{1/2p} \sqrt{n} \big \{ \|\|Z_1-\mathbb{E}(Z_1 \mid \mathcal{F}_0)\|_{\mathbb{H}}\|_{\infty} + \frac{5}{\sqrt{2}} \Delta_q \big \}
\end{equation}
where
\begin{equation*}
\Delta_q=\sum_{j=0}^{q-1} \frac{1}{2^{j/2}} \|\|\mathbb{E}(S_{2^j}\mid \mathcal{F}_0)\|_{\mathbb{H}} \|_{\infty}.
\end{equation*}
\end{prop}
\noindent \textit{Proof of  Proposition \ref{lemhoeffdingstatprocadap}}.
The proof is almost identical to the proof of the corresponding facts in Theorem $1$ of Peligrad, Utev and Wu \cite{peligradutevwu}
if we replace everywhere the absolute value $\mid . \mid $ by $\|.\|_{\mathbb{H}}$, and $C_p$ is here, $C_p=2^{p+1} \Gamma(p+1)$. Also,
we work with the $\mathbf{L}_{\mathbb{H}}^{2p}$-norm. Consequently, we give here, only the crucial inequalities. We note $K=5/\sqrt{2}$. By triangle inequality, notice that
\begin{equation} \label{eqnlemhoeffdingstatprocadap2}
\underset{1 \leq i \leq n}{\max} \|S_i\|_{\mathbb{H}} \leq \underset{1 \leq i \leq n} {\max} \Big\|\sum_{k=1}^i \big(Z_k-\mathbb{E}(Z_k \mid \mathcal{F}_{k-1})\big)\Big  \|_{\mathbb{H}} + \underset{1 \leq i \leq n} {\max} \Big \|\sum_{k=1}^i \mathbb{E}(Z_k \mid \mathcal{F}_{k-1}) \Big  \|_{\mathbb{H}}.
\end{equation}
By the inequality for martingale differences (\ref{eqnlemhoeffdingdiffmart1}), we get
\begin{equation}\label{eqnlemhoeffdingstatprocadap3}
\Big \|\underset{ 1 \leq i \leq n} {\max} \Big \|\sum_{k=1}^i (Z_k-\mathbb{E}(Z_k \mid \mathcal{F}_{k-1})) \Big \|_{\mathbb{H}} \Big \|_{2p} \leq C_p^{1/2p} \sqrt{n} \|\|Z_1-\mathbb{E}(Z_1 \mid \mathcal{F}_0)\|_{\mathbb{H}} \|_{\infty}.
\end{equation}
Moreover, if we start by writing, $n=2m$ or $n=2m+1$, we have
\begin{eqnarray}
\Big \|\underset{1 \leq i \leq n} {\max} \Big \|\sum_{k=1}^i \mathbb{E}(Z_k \mid \mathcal{F}_{k-1}) \Big \|_{\mathbb{H}} \Big \|_{2p} & \leq & \Big \|\underset{1 \leq l \leq m}{\max} \Big \|\sum_{k=1}^{2l} \mathbb{E}(Z_k \mid \mathcal{F}_{k-1}) \Big \|_{\mathbb{H}}\Big\|_{2p} \notag \\
& & + \Big \|\underset{0 \leq l \leq m}{\max} \|\mathbb{E}(Z_{2l+1} \mid \mathcal{F}_{2l}) \|_{\mathbb{H}} \Big \|_{2p}, \label{eqnlemhoeffdingstatprocadap4}
\end{eqnarray}
and
\begin{equation}\label{eqnlemhoeffdingstatprocadap5}
\Big \|\underset{0 \leq l \leq m}{\max}  \|\mathbb{E}(Z_{2l+1} \mid \mathcal{F}_{2l})  \|_{\mathbb{H}} \Big \|_{2p} \leq (m+1)^{1/2p} \|\|\mathbb{E}(Z_1 \mid \mathcal{F}_0)\|_{\mathbb{H}}\|_{2p}.
\end{equation}
For the first term of (\ref{eqnlemhoeffdingstatprocadap4}), we proceed as in the proof of Peligrad, Utev and Wu \cite{peligradutevwu}, to get:
\begin{eqnarray}
\Big \|\underset{1 \leq l \leq m}{\max} \Big \|\sum_{k=1}^{2l} \mathbb{E}(Z_k \mid \mathcal{F}_{k-1}) \Big \|_{\mathbb{H}} \Big \|_{2p} & \leq &  C_p^{1/2p} \sqrt{m} \big \{ 4 \|\|\mathbb{E}(Z_1 \mid \mathcal{F}_0)\|_{\mathbb{H}} \|_{\infty}  \notag \\
& & + K \sqrt{2} \big (\Delta_q-\|\|\mathbb{E}(Z_1 \mid \mathcal{F}_0) \|_{\mathbb{H}} \|_{\infty}\big )\big \}. \label{eqnlemhoeffdingstatprocadap6}
\end{eqnarray}
Consequently, combining (\ref{eqnlemhoeffdingstatprocadap2})-(\ref{eqnlemhoeffdingstatprocadap6}), we obtain the bound (\ref{eqnlemhoeffdingstatprocadap1}).
\hfill $\square$\\ \\
\indent The next proposition is the main tool allowing us to extend Proposition \ref{lemhoeffdingstatprocadap} to non-adapted stationary sequences of $\mathbb{H}$-valued random variables.
\begin{prop} \label{lemhoeffdingstatprocnonadap}
Let $n$, $q$ be integers such that $n \geq 1$, $2^{q-1} \leq n < 2^q$. Assume that $ \|\|Z_0 \|_{\mathbb{H}} \|_{\infty} < \infty$,
and $\mathbb{E}(Z_0 \mid \mathcal{F}_{-1})=0$. Then, for all real $p \geq 1$,
\begin{equation} \label{eqnlemhoeffdingstatprocnonadap1}
\mathbb{E} \Big (\underset{1 \leq i \leq n}{\max} \Big \| \sum_{k=1}^i Z_k \Big \|_{\mathbb{H}}^{2p} \Big) ^{1/2p} \leq (2^{p+1} \Gamma(p+1))^{1/2p} \sqrt{n} \big \{ \|\|\mathbb{E}(Z_0 \mid \mathcal{F}_0)\|_{\mathbb{H}}\|_{\infty} + \frac{2}{\sqrt{2}} \Delta'_q \big \}
\end{equation}
where
\begin{equation*}
\Delta'_q=\sum_{j=0}^{q-1} \frac{1}{2^{j/2}} \|\|S_{2^j}-\mathbb{E}(S_{2^j}\mid \mathcal{F}_{2^j})\|_{\mathbb{H}} \|_{\infty}.
\end{equation*}
\end{prop}
\noindent \textit{Proof of  Proposition \ref{lemhoeffdingstatprocnonadap}}.
Here also, the proof is widely inspired by the proof of Theorem 1 in Peligrad, Utev and Wu \cite{peligradutevwu} and we note always $C_p=2^{p+1} \Gamma(p+1)$. We prove (\ref{eqnlemhoeffdingstatprocnonadap1}) by
induction on $n$. For $n=1$, $q=1$, we clearly have
\begin{equation*}
\|\|Z_1\|_{\mathbb{H}}\|_{2p} \leq \|\|\mathbb{E}(Z_1 \mid \mathcal{F}_1)\|_{\mathbb{H}}\|_{\infty}+ \Delta'_1.
\end{equation*}
Then, assume that the inequality holds for all $n < 2^{q-1}$. Fix $n$ such that, $ 2 ^{q-1} \leq n < 2^q$. By triangle inequality, we obtain that
\begin{equation} \label{eqnlemhoeffdingstatprocnonadap2}
\underset{1 \leq i \leq n}{\max} \Big\|\sum_{k=1}^i Z_k  \Big \|_{\mathbb{H}} \leq \underset{1 \leq i \leq n} {\max} \Big \|\sum_{k=1}^i \mathbb{E}(Z_k \mid \mathcal{F}_k) \Big\|_{\mathbb{H}} + \underset{1 \leq i \leq n}{\max} \Big \|\sum_{k=1}^i \big ( Z_k -\mathbb{E}(Z_k \mid \mathcal{F}_k) \big ) \Big \|_{\mathbb{H}}.
\end{equation}
Since $\mathbb{E}(Z_0 \mid \mathcal{F}_{-1})=0$, we can use the inequality (\ref{eqnlemhoeffdingdiffmart1}) for martingale differences,
\begin{equation}\label{eqnlemhoeffdingstatprocnonadap3}
\Big \|\underset{1 \leq i \leq n}{\max} \Big \|\sum_{k=1}^i \mathbb{E}(Z_k \mid \mathcal{F}_k) \Big \|_{\mathbb{H}}\Big \|_{2p} \leq C_p^{1/2p} \sqrt{n} \|\|\mathbb{E}(Z_0 \mid \mathcal{F}_0) \|_{\mathbb{H}} \|_{\infty}.
\end{equation}
Now, as in Peligrad, Utev and Wu \cite{peligradutevwu}, we write $n=2m$ or $n=2m+1$, for the second term in the right-hand side in (\ref{eqnlemhoeffdingstatprocnonadap2}),
\begin{eqnarray}
\underset{1 \leq i \leq n} {\max}\Big  \| \sum_{k=1}^i \big (Z_k-\mathbb{E}(Z_k \mid \mathcal{F}_k) \big ) \Big \|_{\mathbb{H}} & \leq & \underset{1 \leq l \leq m }{\max} \Big \|\sum_{k=1}^{2l} \big (Z_k-\mathbb{E}(Z_k \mid \mathcal{F}_k) \big ) \Big \|_{\mathbb{H}} \notag \\
& &  + \underset{ 0 \leq l \leq m}{\max} \|Z_{2l+1}-\mathbb{E}(Z_{2l+1} \mid \mathcal{F}_{2l+1}) \|_{\mathbb{H}}, \notag \\
& & \label{eqnlemhoeffdingstatprocnonadap4}
\end{eqnarray}
and
\begin{equation} \label{eqnlemhoeffdingstatprocnonadap5}
\underset{0 \leq l \leq m} {\max} \|Z_{2l+1}-\mathbb{E}(Z_{2l+1} \mid \mathcal{F}_{2l+1})\|_{\mathbb{H}}^{2p} \leq \sum_{l=0}^m  \big \|Z_{2l+1}-\mathbb{E}(Z_{2l+1} \mid \mathcal{F}_{2l+1}) \big \|_{\mathbb{H}}^{2p}.
\end{equation}
For the first term in the right-hand side in (\ref{eqnlemhoeffdingstatprocnonadap4}), we apply the induction hypothesis to
the stationary sequence, $Y_0=Z_0-\mathbb{E}(Z_0 \mid \mathcal{F}_0)+Z_{-1}-\mathbb{E}(Z_{-1} \mid \mathcal{F}_{-1})$, for all $j$ in $\mathbb{Z}$, $Y_j=Y_0 \circ T^{2j}$, the
sigma algebra $\mathcal{G}_0=\mathcal{F}_0$, and the operator $T^2$. Notice that the new filtration becomes $ \{ \mathcal{G}_i : i \in \mathbb{Z} \} $ where $\mathcal{G}_i=\mathcal{F}_{2i}$.
Whence, we have
\begin{equation*}
\Big \|\underset{1 \leq l \leq m}{\max} \Big \|\sum_{k=1}^{2l} \big (Z_k-\mathbb{E}(Z_k \mid \mathcal{F}_k) \big ) \Big \|_{\mathbb{H}} \Big \|_{2p} = \Big \|\underset{1 \leq l \leq m} {\max} \Big \|\sum_{k=1}^l Y_0 \circ T^{2k} \Big \|_{\mathbb{H}} \Big \|_{2p}.
\end{equation*}
Since $ m <2^{q-1}$ and $ \mathbb{E}(Y_0 \mid \mathcal{G}_{-1})=0$, we obtain by the induction hypothesis,
\begin{equation} \label{eqnlemhoeffdingstatprocnonadap6}
\Big\|\underset{1 \leq l \leq m}{\max} \Big\|\sum_{k=1}^l Y _0 \circ T^{2k} \Big \|_{\mathbb{H}} \Big\|_{2p} \leq C_p^{1/2p} \sqrt{m} \big ( \|\|\mathbb{E}(Y_0 \mid \mathcal{G}_0) \|_{\mathbb{H}} \|_{\infty} + \frac{2}{\sqrt{2}} \Delta'_{q-1}(Y) \big ).
\end{equation}
But, $\big \| \big \|\mathbb{E}(Y_0 \mid \mathcal{G}_0) \big \|_{\mathbb{H}} \big \|_{\infty} \leq \|\|Z_0-\mathbb{E}(Z_0 \mid \mathcal{F}_0)\|_{\mathbb{H}} \|_{\infty}$ and rewriting,
\begin{eqnarray}
\Delta'_{q-1}(Y) & = & \sum_{j=0}^{q-2} \frac{1}{2^{j/2}} \Big \| \Big \|\sum_{k=1}^{2^j} Y_k-\mathbb{E}\Big(\sum_{k=1}^{2^j} Y_k \Big | \mathcal{G}_{2^j}\Big) \Big \|_{\mathbb{H}}\Big \|_{\infty} \notag \\
                 & = & \sum_{j=0}^{q-2} \frac{1}{2^{j/2}} \|\|S_{2^{j+1}}-\mathbb{E}(S_{2^{j+1}} \mid \mathcal{F}_{2^{j+1}}) \|_{\mathbb{H}}\|_{\infty} \notag \\
                 & = & \sqrt{2}( \Delta'_q - \|\| Z_0- \mathbb{E}(Z_0 \mid \mathcal{F}_0)\|_{\mathbb{H}} \|_{\infty}), \label{eqnlemhoeffdingstatprocnonadap7}
\end{eqnarray}
we derive that
\begin{eqnarray*}
\Big \|\underset{1 \leq l \leq m}{\max} \Big \|\sum_{k=1}^l Y_0 \circ T^{2k} \Big  \|_{\mathbb{H}} \Big  \|_{2p} &\leq& C_p^{1/2p} \sqrt{m}\big ( \|\|Z_0-\mathbb{E}(Z_0 \mid \mathcal{F}_0)\|_{\mathbb{H}}\|_{\infty} +2 \Delta'_q \notag \\
                                                                                         &     &  -2\|\|Z_0-\mathbb{E}(Z_0\mid \mathcal{F}_0)\|_{\mathbb{H}}\|_{\infty}\big ). \label{eqnlemhoeffdingstatprocnonadap8}
\end{eqnarray*}
Consequently, we conclude that
\begin{eqnarray*}
\Big \|\underset{1 \leq i \leq n}{\max} \Big \|\sum_{k=1}^i Z_k \Big \|_{\mathbb{H}} \Big \|_{2p} & \leq & C_p^{1/2p} \big (\sqrt{n} \|\|\mathbb{E}(Z_0 \mid \mathcal{F}_0) \|_{\mathbb{H}} \|_{\infty} + \sqrt{2m} \frac{2}{\sqrt{2}} \Delta'_q \big ) \notag \\
                                                                              & \leq & C_p^{1/2p} \sqrt{n}\big (\|\|\mathbb{E}(Z_0 \mid \mathcal{F}_0)\|_{\mathbb{H}}\|_{\infty} + \frac{2}{\sqrt{2}} \Delta'_q \big ). \label{eqnlemhoeffdingstatprocnonadap9}
\end{eqnarray*}
\hfill $\square$ \\
\indent Now we give, the main proposition.
\begin{prop} \label{prophoeffding}
Assume that $\|\|X_0\|_{\mathbb{H}}\|_{\infty} < \infty $, then, for all $p \geq 1$,
\begin{eqnarray}
\mathbb{E} \Big (\underset{1 \leq i \leq n}{\max} \|S_i\|_{\mathbb{H}}^{2p} \Big ) & \leq & 2^{p+1} \Gamma(p+1) n^p \Big (\|\|X_0\|_{\mathbb{H}}\|_{\infty} +D \sum_{j=1}^{n} j^{-3/2} \|\|\mathbb{E}(S_j \mid \mathcal{F}_0) \|_{\mathbb{H}} \|_{\infty} \notag \\
                                                                       &      & \quad + D' \sum_{j=1}^{n} j^{-3/2} \|\|S_j - \mathbb{E}(S_j \mid \mathcal{F}_j)\|_{\mathbb{H}}\|_{\infty} \Big )^{2p} \label{eqnprophoeffding}
\end{eqnarray}
where
\begin{equation*}
7 D = 40 \sqrt{2} + 27 \ \ \ \mbox{and} \ \ 7 D'  =  24 \sqrt{2}+12.
\end{equation*}
\end{prop}
\noindent \textit{Proof of  Proposition \ref{prophoeffding}}.
We set $K=\frac{5}{\sqrt{2}}$, $K'=\frac{2}{\sqrt{2}}$ and $C_p=2^{p+1}\Gamma(p+1)$. Let $n$ and $q$ be integers such that $n \geq 1$ and $2^{q-1} \leq n < 2^q$. Let
$$\begin{array}{rclrcl}
\delta_n &=&\sum_{j=1}^{n} j^{-3/2} \|\|\mathbb{E}(S_j \mid \mathcal{F}_0) \|_{\mathbb{H}} \|_{\infty}, & \delta'_n &=&  \sum_{j=1}^n j^{-3/2} \|\|S_j - \mathbb{E}(S_j \mid \mathcal{F}_{j}) \|_{\mathbb{H}} \|_{\infty}, \\
\Delta_q &=&\sum_{j=0}^{q-1} 2^{-j/2} \|\|\mathbb{E}(S_{2^j} \mid \mathcal{F}_0)\|_{\mathbb{H}} \|_{\infty}, & \Delta'_q &=& \sum_{j=0}^{q-1} 2^{-j/2} \|\|S_{2^j}-\mathbb{E}(S_{2^j} \mid \mathcal{F}_{2^j})\|_{\mathbb{H}} \|_{\infty}.
\end{array}$$
We shall prove a slightly stronger inequality,
\begin{equation} \label{eqnprophoeffding1}
\Big \|\mathbb{E} \Big (\underset{1 \leq i \leq n}{\max} \|S_i \|_{\mathbb{H}}\Big ) \Big \|_{2p} \leq C_p^{1/2p} \sqrt{n} (\|\|X_1-\mathbb{E}(X_1 \mid \mathcal{F}_0)\|_{\mathbb{H}}\|_{\infty} +K \Delta_q+K' \Delta'_q).
\end{equation}
Note first that $V_n=\|\| \mathbb{E}(S_n\mid \mathcal{F}_0)\|_{\mathbb{H}}\|_{\infty}$ is a sub-additive sequence as proved by Peligrad and Utev \cite{peligradutev} in Lemma $2.6$ (replace the $\mathbf{L}^2$-norm by the $\mathbf{L}^{\infty}_{\mathbb{H}}$-norm). The sequence $(V_n)_{n \geq 0}$
verifies for all $i,\ j$ in $\mathbb{N}^{\ast}$,
\begin{equation*}
 V_{i+j} \leq V_i +V_j.
\end{equation*}
Whence, using Lemma \ref{appendix} in Appendix with $\tilde{C_1}=\tilde{C_2}=1$, we get
\begin{equation*}
\Delta_q \leq \Big(\frac{4 \sqrt{2}}{7} + \frac{16}{7}\Big) \delta_n.
\end{equation*}
On an other hand, the sequence $V'_n=\|\|S_n -\mathbb{E}(S_n \mid \mathcal{F}_n) \|_{\mathbb{H}}\|_{\infty}$ verifies for all $i,\ j$ in $\mathbb{N}^{\ast}$,
\begin{eqnarray}
V'_{i+j} & \leq & \|\|S_{i+j}- \mathbb{E}(S_{i+j} \mid \mathcal{F}_{i+j})\|_{\mathbb{H}}\|_{\infty} \notag \\
         & \leq & \|\|S_i-\mathbb{E}(S_i \mid \mathcal{F}_i)\|_{\mathbb{H}}\|_{\infty} + \|\|S_{i+j}-S_i-\mathbb{E}(S_{i+j}-S_i \mid \mathcal{F}_{i+j})\|_{\mathbb{H}}\|_{\infty} \notag \\
         &      & \qquad +\|\| \mathbb{E}(S_i \mid \mathcal{F}_i)-\mathbb{E}(S_i \mid \mathcal{F}_{i+j})\|_{\mathbb{H}}\|_{\infty} \notag \\
         & \leq & \|\|S_i-\mathbb{E}(S_i \mid \mathcal{F}_i)\|_{\mathbb{H}}\|_{\infty} + \|\|S_j-\mathbb{E}(S_j \mid \mathcal{F}_j)\|_{\mathbb{H}}\|_{\infty} \notag \\
         &      & \qquad  + \|\|\mathbb{E}(S_i-\mathbb{E}(S_i \mid \mathcal{F}_i)\mid \mathcal{F}_{i+j})\|_{\mathbb{H}}\|_{\infty} \notag \\
         & \leq & 2 V'_i+V'_j. \notag
\end{eqnarray}
Whence, using Lemma \ref{appendix} in Appendix with $\tilde{C_1}=2$ and $\tilde{C_2}=1$, we have
\begin{equation*}
\Delta'_q \leq \Big(\frac{6 \sqrt{2}}{7} + \frac{24}{7}\Big) \delta'_n.
\end{equation*}
Setting  $k_1=\frac{4 \sqrt{2}}{7} + \frac{16}{7}$ and $k_2=\frac{6 \sqrt{2}}{7} + \frac{24}{7}$, we get
\begin{equation*}
\| \|X_1 - \mathbb{E}(X_1 \mid \mathcal{F}_0) \|_{\mathbb{H}} \|_{\infty} + K \Delta_q + K' \Delta'_q \leq \|\|X_1 \|_{\mathbb{H}}\|_{\infty} + ( K k_1 +1 ) \delta_n +K' k_2 \delta'_n.
\end{equation*}
Since $(\ref{eqnprophoeffding1})$ implies $(\ref{eqnprophoeffding})$, it remains to prove $(\ref{eqnprophoeffding1})$.\\
\indent By triangle inequality, we obtain that
\begin{eqnarray}
 \Big\|\underset{1 \leq i \leq n } {\max} \Big\|\sum_{k=1}^i X_k \Big\|_{\mathbb{H}} \Big\|_{2p}
& \leq & \Big\|\underset{1 \leq i \leq n} {\max} \Big\|\sum_{k=1}^i \big(X_k -\mathbb{E}(X_k \mid \mathcal{F}_k)\big) \Big\|_{\mathbb{H}} \Big\|_{2p} \notag \\
&& \quad \quad  + \Big\|\underset{1 \leq i \leq n}{\max} \Big\|\sum_{k=1}^i \mathbb{E}(X_k \mid \mathcal{F}_k ) \Big\|_{\mathbb{H}} \Big\|_{2p}.  \label{eqnprophoeffding2}
\end{eqnarray}
Applying Proposition \ref{lemhoeffdingstatprocadap}, we derive
\begin{equation}\label{eqnprophoeffding4}
\Big \|\underset{1 \leq i \leq n}{\max} \Big \|\sum_{k=1}^i \mathbb{E}(X_k \mid \mathcal{F}_k) \Big\|_{\mathbb{H}} \Big\|_{2p} \leq C_p ^{1/2p} \sqrt{n}\big(\|\|\mathbb{E}(X_1 \mid \mathcal{F}_1)-\mathbb{E}(X_1 \mid \mathcal{F}_0) \|_{\mathbb{H}}\|_{\infty} +K \Delta^{\ast}_q \big )
\end{equation}
where
\begin{equation*}
\Delta^{\ast}_q =  \sum_{j=0}^{q-1} 2^{-j/2} \Big \| \Big \|\mathbb{E}\Big(\sum_{k=1}^{2^j} \mathbb{E}(X_k \mid \mathcal{F}_k) \Big | \mathcal{F}_0 \Big)\Big  \|_{\mathbb{H}} \Big\|_{\infty}=\Delta_q.
\end{equation*}
On the other hand, Proposition \ref{lemhoeffdingstatprocnonadap} gives
\begin{equation}\label{eqnprophoeffding5}
\Big \|\underset{1\leq i \leq n} {\max} \Big \|\sum_{k=1}^i \big ( X_k -\mathbb{E}(X_k \mid \mathcal{F}_k) \big ) \Big \|_{\mathbb{H}} \Big \|_{2p} \leq C_p^{1/2p} \sqrt{n} \{ \|\|\mathbb{E}(X_0-\mathbb{E}(X_0 \mid \mathcal{F}_0)\mid \mathcal{F}_0)\|_{\mathbb{H}} \|_{\infty} +K' \Delta'^{\ast}_q \}
\end{equation}
where
\begin{equation*}
\Delta'^{\ast}_q  =  \sum_{j=0}^{q-1} 2^{-j/2} \Big \| \Big \|\sum_{k=1}^{2^j} \big ( \{X_k-\mathbb{E}(X_k \mid \mathcal{F}_k)\}-\mathbb{E}(X_k - \mathbb{E}(X_k \mid \mathcal{F}_k) \mid \mathcal{F}_{2^j}) \big ) \Big \|_{\mathbb{H}} \Big \|_{\infty}=\Delta'_q.
\end{equation*}
Combining $(\ref{eqnprophoeffding4})$ and $(\ref{eqnprophoeffding5})$ in $(\ref{eqnprophoeffding2})$, $(\ref{eqnprophoeffding1})$ follows.
\hfill $\square$
\subsubsection{Proof of Theorem \ref{theohoeffding}}
$\newline$
Now, assume that $p$ is an integer, and $p \geq 1$.\\
Let
\begin{equation*}
B =  \|\|X_0 \|_{\mathbb{H}} \|_{\infty} + \Big (D \sum_{j=1}^{n} \frac{1}{j^{3/2}} \| \| \mathbb{E}(S_j \mid \mathcal{F}_0) \|_{\mathbb{H}} \|_{\infty}
+D' \sum_{j=1}^{n} \frac{1}{j^{3/2}} \| \|S_j - \mathbb{E}(S_j\mid \mathcal{F}_j) \|_{\mathbb{H}}\|_{\infty} \Big )
\end{equation*}
for some constants $D>0$ and $D'>0$ defined in Proposition \ref{prophoeffding}.
We can use the approach of the proof of Theorem $2.4$ in Rio \cite{rio}, because
\begin{eqnarray*}
\mathbb{E}\Big (\underset{1 \leq i \leq n}{\max} \|S_i\|^{2p}_{\mathbb{H}} \Big ) & \leq & 2^{p+1} p! B^{2p} n^p \\
                                                                       & \leq & 2 (2p-1)!! (2nB^2)^p.
\end{eqnarray*}
Consequently , if we use the notation of the proof of Theorem $2.4$ in Rio \cite{rio}, the constant $A$ is here,
\begin{equation*}
A=\frac{x^2}{4 n B^2},
\end{equation*}
and with the estimation given in Rio \cite{rio} (page $42$),
\begin{equation*}
\mathbb{P} \Big (\underset{1 \leq i \leq n}{\max} \|S_i\|_{\mathbb{H}} \geq x \Big ) \leq 2 \sqrt{e} \exp(-A).
\end{equation*}
Taking $C=\max \{D,D'\}$, we obtain exactly Theorem \ref{theohoeffding}.
\hfill $\square$
\subsection{MDP for martingale differences}
$\newline$
\indent Our main proposition is a generalization of a result of Theorem $3.1$ in Puhalskii \cite{puhalskii} to $\mathbb{H}$-valued random variables.
\begin{prop} \label{propMDPdiffmart}
Let $a_n$ be a sequence of positive numbers satisfying $a_n \rightarrow 0$ and $na_n \underset{n \rightarrow \infty}{\rightarrow} + \infty$.
Let $k_n$ be an increasing sequence of integers going to infinity and $ \{ d_{j,n} \}_{1 \leq j \leq k_n} $ be a triangular array of martingale differences, with values in $\mathbb{H}$, such that
\begin{equation} \label{eqnMDPdiffmart1}
\forall \  1 \leq j \leq k_n, \ \| d_{j,n} \| _{\mathbf{L}^{\infty}_{\mathbb{H}}} \leq \beta_n \sqrt{n a_n} \ \ \mbox{with $\underset{n \longrightarrow \infty}{\beta_n \longrightarrow 0}$}.
\end{equation}
Assume that, there exists $Q \in \mathcal{S}(\mathbb{H})$ such that:
\begin{enumerate}
\item[$i.$] for all $k$, $l$ in $\mathbb{N}^{\ast}$ and $\delta>0$,
\begin{equation} \label{eqnMDPdiffmart2}
\underset{n\longrightarrow \infty}{\limsup} \ a_n \ \log \mathbb{P} \Big ( \Big | \frac{1}{n} \sum_{j=1}^{k_n} \mathbb{E}(<d_{j,n},e_k>_{\mathbb{H}} <d_{j,n},e_l>_{\mathbb{H}} \mid \mathcal{F}_{j-1,n})-<Qe_k,e_l>_{\mathbb{H}} \Big | > \delta \Big )=-\infty,
                        \end{equation}
\item[$ii.$] for all $\delta >0$,
\begin{equation} \label{eqnMDPdiffmart3}
 \underset{n\longrightarrow \infty}{\limsup} \ a_n \ \log \mathbb{P} \Big ( \Big | \frac{1}{n} \sum_{j=1}^{k_n} \mathbb{E}(\| d_{j,n}\| ^2_{\mathbb{H}} \mid \mathcal{F}_{j-1,n})-\mbox{Tr}(Q) \Big | > \delta  \Big )=-\infty.
\end{equation}
\end{enumerate}
Then $\{W_n(t)=\frac{1}{\sqrt{n}} \sum_{j=1}^{[k_nt]} d_{j,n}+\frac{1}{\sqrt{n}}(k_nt-[k_nt])d_{[k_nt]+1,n} : t \in [0,1] \}$ satisfies the MDP in $C_{\mathbb{H}}([0,1])$ with speed $a_n$ and the good rate function,
\begin{equation} \label{eqnMDPdiffmart4}
I(\phi) = \left  \lbrace \begin{array}{l} \int_0^1 \Lambda ^{\ast}(\phi'(t))\, dt \ \mbox{if $\phi \in \mathcal{AC}_0([0,1])$} \\ \\
                    \infty \ \mbox{otherwise}
\end{array}
\right.
\end{equation}
where $\Lambda ^{\ast}$ is defined by
\begin{equation} \label{eqnMDPdiffmart5}
 \Lambda ^{\ast}(x)=\underset{y \in {\mathbb{H}}}{\sup} \big (<y,x>_{\mathbb{H}}-\frac{1}{2}<y,Qy>_{\mathbb{H}} \big).
\end{equation}
\end{prop}
\noindent \textit{Proof of  Proposition \ref{propMDPdiffmart}}.
$\newline$
Firstly, we need some notations.
\begin{nota} \label{Notation3}
For all integer $m \geq 1$, let $P^m$ be the projection on the first $m$
components of the orthonormal basis, $(e_i)_{1 \leq i \leq m}$, in $ \mathbb{H} $ then
\begin{equation*} \label{eqnnotation3}
d_{j,n}^m=P^m(d_{j,n}), \ r_{j,n}^m=(I-P^m)d_{j,n}.
\end{equation*}
where $I$ is the identity operator.
\end{nota}
\indent Let $ \{ d_{j,n} \}_{1\leq j \leq k_n}$ be a $\mathbb{H}$-valued triangular array of martingale differences. We start by proving that $ \{d_{j,n}^m\}_{1 \leq j \leq k_n}$, which is
a $\mathbb{R}^m$-valued triangular array of martingale differences satisfies the conditions of Theorem $3.1$ of Pulhalskii \cite{puhalskii} (see also Djellout \cite{djellout}, Proposition $1$). \\
\indent The conditions $(\ref{eqnMDPdiffmart1})$ and $(\ref{eqnMDPdiffmart2})$ imply conditions $i)$ and $ii)$ of Proposition 1 in Djellout \cite{djellout}. \\
\indent Consequently, $\{W_{n}^{m}(t)=\frac{1}{\sqrt{n}} \sum_{j=1}^{[k_nt]} d_{j,n}^m+\frac{1}{\sqrt{n}}(k_nt-[k_nt])d^m_{[k_nt]+1,n}: t \in[0,1]\} $
satisfies  the MDP, with the good rate function, $I_m(.)$,
\begin{equation*} \label{eqnMDPdiffmart8}
I_m(\phi) = \left  \lbrace \begin{array}{l} \int_0^1 \Lambda_m ^{\ast}(\phi'(t)) \, dt \ \mbox{if $\phi \in \mathcal{AC}_0([0,1])$} \\ \\
                    \infty \ \mbox{otherwise},
\end{array}
\right.
\end{equation*}
where $\Lambda_m ^{\ast}$ is:
\begin{equation*} \label{eqnMDPdiffmart9}
\forall \ x \in \mathbb{H}, \  \Lambda_m ^{\ast}(x)=\underset{y \in {\mathbb{H}}}{\sup} \big (<P^m y,P^m x>_{\mathbb{H}}-\frac{1}{2}<P^m y,Q P^m y>_{\mathbb{H}} \big ).
\end{equation*}
\indent By using Theorem $4.2.13$ in Dembo and Zeitouni \cite{dembozeitouni}, it remains to prove, that for any $\eta >0$,
\begin{equation*} \label{eqn394}
 \underset{m \longrightarrow \infty}{\limsup} \  \underset{n \longrightarrow \infty}{\limsup} \ a_n \ \log \ \mathbb{P} \Big ( \underset{1 \leq j \leq k_n} {\max} \sqrt{\frac{a_n}{n}} \Big  \| \sum_{k=1}^{j} r_{k,n}^m \Big \|_{\mathbb{H}} > \eta \Big )=-\infty.
\end{equation*}
Notice that, for all $\eta>0$,
\begin{multline*} \label{eqn395}
a_n \ \log \ \mathbb{P} \Big ( \underset{1 \leq j \leq k_n}{\max} \sqrt{\frac{a_n}{n}} \Big \|  \sum_{k=1}^j r_{k,n}^m \Big \|_{\mathbb{H}} > \eta \Big )  \\
\leq  a_n \ \log \ \Big ( \mathbb{P} \Big ( \Big \{ \underset{1 \leq j \leq k_n} {\max} \sqrt{\frac{a_n}{n}} \Big \| \sum_{k=1}^j r_{k,n}^m  \Big \|_{\mathbb{H}} > \eta  \Big \} \cap \Big \{ \Big | \frac{1}{n} \sum_{k=1}^{k_n} \mathbb{E}( \|r_{k,n}^m \|_{\mathbb{H}}^2 \mid \mathcal{F}_{k-1,n}) \\
 \qquad  -\sum_{p=m+1}^{\infty} < Qe_p,e_p>_{\mathbb{H}} \Big | \leq \varepsilon \Big \} \Big )
  + \mathbb{P} \Big ( \Big | \frac{1}{n} \sum_{k=1}^{k_n} \mathbb{E}( \|r_{k,n}^m \|_{\mathbb{H}}^2 \mid \mathcal{F}_{k-1,n}) - \sum_{p=m+1}^{\infty} < Qe_p,e_p>_{\mathbb{H}} \Big | > \varepsilon \Big )\Big),
\end{multline*}
where $\varepsilon>0$. \\
With the notations
\begin{multline*}
A(n,m, \eta , \varepsilon):=\mathbb{P} \Big ( \Big \{ \underset{ 1 \leq j \leq k_n}{\max} \sqrt{ \frac{a_n}{n}} \Big \|  \sum_{k=1}^j r_{k,n}^m  \Big \| _{\mathbb{H}} > \eta \Big \} \\
 \cap \Big \{  \Big | \frac{1}{n} \sum_{k=1}^ {k_n}  \mathbb{E} ( \| r_{k,n}^m \|^2 _{\mathbb{H}} \mid \mathcal{F}_{k-1,n} ) - \sum_{p=m+1}^{\infty} < Q \ e_p, e_p>_{\mathbb{H}} \Big | \leq \varepsilon \Big \}  \Big),
\end{multline*}
and
\begin{equation*}
B(n,m, \varepsilon):= \mathbb{P} \Big(  \Big | \frac{1}{n} \sum_{j=1}^{k_n}
\mathbb{E} ( \| r_{j,n}^m \| ^2_{\mathbb{H}} \mid \mathcal{F}_{j-1,n} ) -
\sum_{p=m+1}^{\infty} < Q \ e_p, e_p>_{\mathbb{H}} \Big | > \varepsilon \Big),
\end{equation*}
we derive
\begin{equation*}
 a_n  \mbox{log} \ \mathbb{P} \Big ( \underset{1 \leq j \leq k_n}{\max} \sqrt{\frac{a_n}{n}} \Big \|  \sum_{k=1}^j r_{k,n}^m  \Big \| _{\mathbb{H}} > \eta  \Big)
  \leq  a_n  \log \ \big \{ A(n, m, \eta, \varepsilon) + B(n, m, \varepsilon) \big \}. \ \
\end{equation*}
Now notice
\begin{equation*} \begin{array}{rcl} \label{3008}
& & a_n  \log ( B(n, m, \varepsilon) ) \\ \\
& \leq & a_n  \log \ \big \{ \mathbb{P} \big (  \big | \frac{1}{n} \sum_{j=1}^{k_n}  \mathbb{E} ( \| d_{j,n} \| ^2 _{\mathbb{H}} \mid \mathcal {F}_{j-1,n} ) - \mbox{Tr}(Q) \big | >\frac{\varepsilon }{2} \big) \\ \\
&      &   + \mathbb{P} \big ( \big  | \frac{1}{n} \sum_{j=1}^{k_n}
\mathbb{E} ( \| d^m_{j,n} \| ^2_{\mathbb{H}} \mid \mathcal{F}_{j-1,n} ) -
\sum_{p=1}^{m} < Q \ e_p, e_p>_{\mathbb{H}} \big | > \frac{\varepsilon }{2} \big) \big \}. \
\
\end{array}
\end{equation*}
Using (\ref{eqnMDPdiffmart2}) and (\ref{eqnMDPdiffmart3}), it follows
\begin{equation*} \label{eqn308}
 \underset{m \longrightarrow  \infty} { \limsup} \ \underset{ n \longrightarrow  \infty} { \limsup} \ a_n \ \mbox{log} (B(n, m, \varepsilon))=- \infty. \ \
\end{equation*}
With the notations
\begin{equation*}
 C(n,m, \eta) := \Big \{ \underset{ 1 \leq j \leq k_n}{\max} \sqrt{ \frac{a_n}{n}} \Big \| \sum_{k=1}^j r_{k,n}^{m}  \Big \| _{\mathbb{H}} > \eta \Big \}, \ \
\end{equation*}
and
\begin{equation*}
 D(n,m, \varepsilon) := \Big \{  \Big | \frac{1}{n}  \sum_{j=1}^ {k_n}  \mathbb{E} ( \| r_{j,n}^m \|^2 _{\mathbb{H}} \mid \mathcal{F}_{j-1,n} )- \sum_{p=m+1}^{\infty} < Q \ e_p, e_p>_{\mathbb{H}} \Big | \leq \varepsilon \Big\}, \ \
\end{equation*}
applying Theorem $5.1$ (inequality $(5.2)$), in Kallenberg and Sztencel \cite{kallenbergsztencel91} or Theorem $3.4$ in Pinelis \cite{pinelis94} to the martingale difference,
$$U_i= r_{i,n}^m \mathbf{1}_{\{\sum_{j=1}^i \mathbb{E}(\|r_{j,n}^m\|_{\mathbb{H}} \mid \mathcal{F}_{j-1,n}) \leq n(\varepsilon+\sum_{p=m+1}^{\infty} <Qe_p,e_p>_{\mathbb{H}})\}},$$ we obtain that
\begin{multline*}
 a_n \ \log \ \mathbb{P}( C(n,m, \eta) \cap D(n,m, \varepsilon) )\\
\begin{split}
& \leq a_n \log(2)+ \ a_n  \log  \Big( \exp \Big( - \frac{ \eta ^2 n }{ 2 a_n \big (n \varepsilon  + n \sum_{p=m+1}^{\infty} < Q \ e_p, e_p>_{\mathbb{H}} \big) + \frac{2}{3} a_n \beta _n \sqrt{n a_n } \sqrt{n} \frac{\eta }{\sqrt{a_n}}} \Big) \Big) \ \ \  \\
&  \leq a_n \log(2)   - \frac{\eta ^2}{ 2\varepsilon + 2\sum_{p=m+1}^{\infty} < Q \
e_p, e_p>_{\mathbb{H}}+ \frac{2}{3} \beta_n  \eta  }. \ \
\end{split}
\end{multline*}
Since $Q$ has a finite trace, it follows
\begin{equation*}
 \underset{m\rightarrow \infty}{\lim} \ \underset{n \rightarrow \infty}{\lim} \ \underset{ \varepsilon \rightarrow 0}{\lim} - \frac{\eta ^2}{2 \varepsilon + 2\sum_{p=m+1}^{\infty} < Q \ e_p, e_p>_{\mathbb{H}} + \frac{2}{3} \beta_n  \eta }=- \infty.  \ \
\end{equation*}
Consequently, we conclude by Theorem $4.2.13$, in Dembo and Zeitouni \cite{dembozeitouni} that
$ \{n^{-1/2} \sum_{j=1}^{[k_nt]} d_{j,n}+\frac{1}{\sqrt{n}}(k_nt-[k_nt])d_{[k_nt]+1,n}:  t \in[0,1]\} $
satisfies the MDP in $C_{\mathbb{H}}([0,1])$.
The rate function is the same that the i.i.d gaussian random variable with mean $0$ and covariance $Q$, therefore equal to
\begin{equation*} \label{eqn315}
\forall \ x \in \mathbb{H}, \  \Lambda^{\ast}(x) = \underset{y \in \mathbb{H}}{\sup} \big( <y,x>_{\mathbb{H}} - \frac{1}{2} <y,Q y>_{\mathbb{H}} \big). \ \
\end{equation*}
\hfill $\square$

\subsection{MDP for stationary sequences}
\subsubsection{Proof of Theorem \ref{theomdpfonctionnel}}
$\newline$
\indent The proof of Theorem \ref{theomdpfonctionnel} uses the same arguments as in the proof of Theorem $1$ in
Dedecker, Merlev\`ede, Peligrad and Utev \cite{dedeckermerlevedepeligradutev}, but for a $\mathbb{H}$-valued non-adapted sequences.\\
\indent Let $m_n=o(\sqrt{na_n})$,
and $k_n=[n/m_n]$ (where, as before, $[x]$ denotes the integer part of $x$). \\
\indent We divide the variables in blocks of size $m_n$ and make the sums in each block,
\begin{equation*}
X_{i,m_n}=\sum_{j=(i-1)m_n+1}^{im_n} X_j, \ i \geq 1.
\end{equation*}
Then, we construct the martingales,
\begin{eqnarray*}
M_{k_n}^{(m_n)} & = &  \sum_{i=1}^{k_n} \big ( \mathbb{E}(X_{i,m_n} \mid \mathcal{F}_{im_n})-\mathbb{E}(X_{i,m_n} \mid \mathcal{F}_{(i-1)m_n}) \big )   \\
          & := & \sum_{i=1}^{k_n} D_{i,m_n},
\end{eqnarray*}
and we define the process $\{ M_{k_n}^{(m_n)}(t):  t \in [0,1] \}$ by
\begin{equation*}
M_{k_n}^{(m_n)}(t) :=M_{[k_n t]}^{(m_n)}+\frac{1}{\sqrt{n}}(k_nt-[k_nt])D_{[k_nt]+1,m_n}.
\end{equation*}
Now, we shall use Proposition \ref{propMDPdiffmart}, applied with $d_{j,n}=D_{j,m_n}$, and verify the conditions (\ref{eqnMDPdiffmart2}) and (\ref{eqnMDPdiffmart3}). \\
\indent We start by proving (\ref{eqnMDPdiffmart2}). By stationarity, it is enough to prove that, for all $k$, $l \geq 1$,
\begin{equation} \label{eqnfonct4}
\underset{n \longrightarrow \infty}{\limsup} \  \Big \| \frac{1}{m_n} \mathbb{E}(<D_{1,m_n},e_k>_{\mathbb{H}} <D_{1,m_n},e_l>_{\mathbb{H}} \mid \mathcal{F}_0)-<Q e_k, e_l>_{\mathbb{H}} \Big \|_{\infty} =0.
\end{equation}
But, we notice
\nolinebreak
\begin{eqnarray*}
& & \mathbb{E}(<D_{1,m_n},e_k> _{\mathbb{H}} < D_{1,m_n},e_l>_{\mathbb{H}} \mid \mathcal{F}_0) \notag \\
& = &  \mathbb{E}(< \mathbb{E}(X_{1,m_n} \mid \mathcal{F}_{m_n}), e_k>_{\mathbb{H}} < \mathbb{E}(X_{1,m_n} \mid \mathcal{F}_{m_n}),e_l>_{\mathbb{H}} \mid \mathcal{F}_{0}) \notag \\
& & \quad - < \mathbb{E}(X_{1,m_n} \mid \mathcal{F}_{0}), e_k>_{\mathbb{H}} < \mathbb{E}(X_{1,m_n} \mid \mathcal{F}_{0}), e_l>_{\mathbb{H}}, \notag
\end{eqnarray*}
\nolinebreak
thus
\nolinebreak
\begin{align*}
& \| \mathbb{E}(<D_{1,m_n},e_k>_{\mathbb{H}} <D_{1,m_n},e_l>_{\mathbb{H}} \mid \mathcal{F}_{0}) - <Qe_k,e_l>_{\mathbb{H}} \|_{\infty} \notag & \\
& \leq \| \mathbb{E}(<\mathbb{E}(S_{m_n} \mid \mathcal{F}_{m_n}),e_k>_{\mathbb{H}} < \mathbb{E}(S_{m_n} \mid \mathcal{F}_{m_n}),e_l>_{\mathbb{H}} \mid \mathcal{F}_{0})-<Qe_k,e_l>_{\mathbb{H}} \|_{\infty} \notag & \\
& \qquad +\| < \mathbb{E}(S_{m_n} \mid \mathcal{F}_{0}),e_k>_{\mathbb{H}} < \mathbb{E}(S_{m_n} \mid \mathcal{F}_{0}),e_l>_{\mathbb{H}} \|_{\infty}. \notag &
\end{align*}
\nolinebreak
By triangle inequality and Cauchy-Schwarz inequality, we have
\begin{align*}
&  \Big \| \frac{1}{m_n} \mathbb{E}(< \mathbb{E}(S_{m_n} \mid \mathcal{F}_{m_n}),e_k>_{\mathbb{H}} < \mathbb{E}(S_{m_n} \mid \mathcal{F}_{m_n}),e_l>_{\mathbb{H}} \mid \mathcal{F}_{0})-<Qe_k,e_l>_{\mathbb{H}} \Big \|_{\infty} \notag & \\
& \leq \frac{1}{m_n} \|\|\mathbb{E}(S_{m_n} \mid \mathcal{F}_{m_n})-S_{m_n} \|_{\mathbb{H}}\|_{\infty}^2 \notag & \\
& \quad + \frac{2}{m_n} \|\|\mathbb{E}(S_{m_n} \mid \mathcal{F}_{m_n})-S_{m_n}\|_{\mathbb{H}}\|_{\infty} \|\sqrt{\mathbb{E}(\|S_{m_n}\|_{\mathbb{H}}^2 \mid \mathcal{F}_0)}\|_{\infty} \notag & \\
& \quad + \Big \| \frac{1}{m_n} \mathbb{E}(<S_{m_n},e_k>_{\mathbb{H}}<S_{m_n},e_l>_{\mathbb{H}} \mid \mathcal{F}_0)-<Qe_k,e_l>_{\mathbb{H}} \Big\|_{\infty}. \notag &
\end{align*}
By using Lemma \ref{appendix} in Appendix, and the hypothesis (\ref{eqnfonct1}),
we deduce (\ref{eqnfonct4}).\\
\indent Now, to prove (\ref{eqnMDPdiffmart3}), by stationarity, we have to verify
\begin{equation} \label{eqnfonct6}
\underset{n \longrightarrow \infty}{\limsup} \Big \|\frac{1}{m_n} \mathbb{E}(\|D_{1,m_n}\|_{\mathbb{H}}^2 \mid \mathcal{F}_{0})- \mbox{Tr}(Q) \Big \|_{\infty} =0.
\end{equation}
Notice that
\begin{equation*}
 \mathbb{E}(\|D_{1,m_n}\|^2_{\mathbb{H}} \mid \mathcal{F}_{0})
\leq \mathbb{E}(\|S_{m_n} \|_{\mathbb{H}}^2 \mid \mathcal{F}_{0})- \|\mathbb{E}(S_{m_n} \mid \mathcal{F}_{0})\|_{\mathbb{H}}^2,
\end{equation*}
thus
\begin{equation*}
\Big \| \frac{1}{m_n}  \mathbb{E}( \|D_{1,m_n} \|_{\mathbb{H}}^2 \mid \mathcal{F}_{0})- \mbox{Tr}(Q) \Big \|_{\infty} \leq \Big \|\frac{1}{m_n}  \mathbb{E}( \|S_{m_n} \|^2_{\mathbb{H}} \mid \mathcal{F}_{0}) - \mbox{Tr}(Q) \Big  \|_{\infty} + \frac{1}{m_n} \| \|\mathbb{E}(S_{m_n} \mid \mathcal{F}_{0})\|_{\mathbb{H}}^2 \|_{\infty}.
\end{equation*}
By using Lemma \ref{appendix} in Appendix, and the hypothesis (\ref{eqnfonct2}),
 we deduce (\ref{eqnfonct6}).\\
\indent To finish the proof, it remains to prove,  that for all $\delta > 0$,
$\nobreak$
\begin{equation} \label{eqnfonct7}
\underset{n \longrightarrow \infty}{\limsup} \ a_n \log \mathbb{P} \Big (\sqrt{\frac{a_n}{n}} \underset{ t \in [0,1]}{\sup} \|S_{[nt]}-M_{[k_nt]}^{(m_n)} \|_{\mathbb{H}}  \geq \delta \Big )=-\infty
\end{equation}
and
\begin{equation} \label{eqnfonct71}
\underset{n \longrightarrow \infty}{\limsup} \ a_n \log \mathbb{P} \Big (\sqrt{\frac{a_n}{n}} \underset{ t \in [0,1]}{\sup} \|(k_nt-[k_nt])D_{[k_nt]+1,m_n}-(nt-[nt])X_{[nt]+1} \|_{\mathbb{H}}  \geq \delta \Big )=-\infty.
\end{equation}
$(\ref{eqnfonct71})$ holds since $m_n=o(\sqrt{a_nn})$ and the random variables are bounded. We turn now to the proof of $(\ref{eqnfonct7})$.
Notice that
$\nobreak$
\begin{eqnarray}
\underset{t \in [0,1]}{\sup} \|S_{[nt]}-M_{[k_nt]}^{(m_n)} \|_{\mathbb{H}} & \leq & \underset{t \in [0,1]}{\sup} \Big \|\sum_{i=[k_n t]m_n+1}^{[nt]} X_i \Big \|_{\mathbb{H}} + \underset{ t \in [0,1]}{\sup} \Big \|\sum_{i=1}^{[k_n t]} \big ( X_{i,m_n}-\mathbb{E}(X_{i,m_n} \mid \mathcal{F}_{im_n}) \big ) \Big\|_{\mathbb{H}} \notag \\
& & \quad + \underset{t \in [0,1]}{\sup} \Big \| \sum_{i=1}^{[k_nt]} \mathbb{E}(X_{i,m_n} \mid \mathcal{F}_{(i-1)m_n}) \Big\|_{\mathbb{H}} \notag \\
                     & \leq & o(\sqrt{na_n})+\underset{1 \leq j \leq k_n}{\max} \Big \|\sum_{i=1}^{j} \big ( X_{i,m_n}-\mathbb{E}(X_{i,m_n} \mid \mathcal{F}_{im_n}) \big ) \Big \|_{\mathbb{H}} \notag \\
                     & & \quad + \underset{1 \leq j \leq k_n}{\max} \Big \| \sum_{i=1}^{j} \mathbb{E}(X_{i,m_n} \mid \mathcal{F}_{(i-1)m_n}) \Big \|_{\mathbb{H}}. \label{eqnfonct00}
\end{eqnarray}
For the last term of the right-hand side in (\ref{eqnfonct00}), we use the same arguments as in the proof of Theorem 1 in Dedecker, Merlev\`ede, Peligrad and Utev \cite{dedeckermerlevedepeligradutev}, so we give
only the proof of the non-adapted term, i.e. the second term of the right-hand side of inequality (\ref{eqnfonct00}). \\
\indent We apply Theorem \ref{theohoeffding} to the stationary sequence, $Y_{0,m_n}=X_{0,m_n}-\mathbb{E}(X_{0,m_n} \mid \mathcal{F}_{0})$,
and $Y_{i,m_n}=Y_{0,m_n}\circ T^{im_n}$. Notice that the new filtration becomes $\{ \mathcal{G}_i,i \in \mathbb{Z} \} $ where $\mathcal{G}_0=\mathcal{F}_0$, and $\mathcal{G}_i =T^{-(im_n)}(\mathcal{G}_0)$. \\
Consequently, we have
\begin{align}
& a_n \log \mathbb{P} \Big (\sqrt{\frac{a_n}{n}} \underset{1 \leq j \leq k_n}{\max}  \Big \|\sum_{i=1}^{j} \big ( X_{i,m_n}-\mathbb{E}(X_{i,m_n} \mid \mathcal{F}_{im_n}) \big ) \Big \|_{\mathbb{H}} \geq \delta \Big ) \notag & \\
& = a_n \log \mathbb{P} \Big (\sqrt{\frac{a_n}{n}} \underset{1 \leq j \leq k_n}{\max} \Big \|\sum_{i=1}^{j} Y_{i,m_n} \Big \|_{\mathbb{H}} \geq \delta \Big ) \notag & \\
& \leq a_n \log(2 \sqrt{e})-\frac{ \delta^2}{ 4  \frac{1}{m_n} E(n,\delta)^2 }, \label{eqnfonct9} &
\end{align}
where
\begin{equation*}
E(n,\delta):=\|\|S_{m_n}-\mathbb{E}(S_{m_n} \mid \mathcal{F}_{m_n}) \|_{\mathbb{H}}\|_{\infty} + C \sum_{j=1}^{\infty} \frac{1}{j^{3/2}} \|\|S_{jm_n}-\mathbb{E}(S_{jm_n} \mid \mathcal{F}_{jm_n}) \|_{\mathbb{H}}\|_{\infty},
\end{equation*}
with $C$ is the positive constant defined in Theorem \ref{theohoeffding}. \\
\indent We conclude by using Lemma \ref{appendix} in Appendix, and  the inequality (\ref{eqnfonct9}) converges to $0$, when $n \longrightarrow \infty$.
\hfill $\square$
\subsection{Proof of Corollary \ref{corollaire}}
$\newline$
\indent The proof of Corollary \ref{corollaire} uses the same arguments as in the proof of Corollary $2$ in Dedecker, Merlev\`ede, Peligrad and Utev \cite{dedeckermerlevedepeligradutev}
but for a non-adapted stationary $\mathbb{H}$-valued sequence. \\
\indent By triangle inequality and changing the order of summation, (\ref{conditionfonctionnelcor}) implies (\ref{conditionfonctionnel}).
\subsubsection{A technical lemma}
\begin{lem}\label{lemtech}
Assume that $ \|\|X_0\|_{\mathbb{H}} \|_{\infty} < \infty $.
Let $n$ be a diadic integer, $n=2^q$. Then
\begin{eqnarray}
\|\mathbb{E}(\|S_n\|_{\mathbb{H}} ^2 \mid \mathcal{F}_0) \|_{\infty} & \leq &  n ( \|\mathbb{E}(\|X_1\|_{\mathbb{H}}^2 \mid \mathcal{F}_0)\|_{\infty} + \frac{1}{2} \Delta_q +\frac{1}{2} \Delta'_q)^2 \label{lemtech1} \\
                                                                     & \leq & n \Delta_{\infty}^2 \notag
\end{eqnarray}
where $\Delta_q$, $\Delta'_q$ are, respectively, defined as in Proposition \ref{lemhoeffdingstatprocadap} and in Proposition \ref{lemhoeffdingstatprocnonadap} and
\begin{equation*}
\Delta_{\infty}=\|\mathbb{E}(\|X_1\|_{\mathbb{H}}^2 \mid \mathcal{F}_0)\|_{\infty} +\frac{1}{2} \Delta_q+\frac{1}{2} \Delta'_q.
\end{equation*}
\end{lem}
\noindent \textit{Proof of Lemma \ref{lemtech}}.
As in the proof of Proposition $2.1$ in Peligrad and Utev \cite{peligradutev}, we prove Lemma  \ref{lemtech} by induction on $q$. \\
\indent Obviously, (\ref{lemtech1}) is true for $q=0$. Assume now, that (\ref{lemtech1}) holds for all diadic integers $n \leq 2^{q-1}$. \\
\indent Writing $S_{2^q}= S_{2^{q-1}}+S_{2^q}-S_{2^{q-1}}$, notice that
\begin{equation*}
\|S_{2^q}\|^2_{\mathbb{H}}  = \|S_{2^{q-1}} \|_{\mathbb{H}}^2 + \|S_{2^q}-S_{2^{q-1}}\|_{\mathbb{H}}^2 + 2 <S_{2^{q-1}},S_{2^q}-S_{2^{q-1}}>_{\mathbb{H}}.
\end{equation*}
By stationarity, we have
\begin{align}
& \|\mathbb{E}(\|S_{2^q}\|_{\mathbb{H}}^2 \mid \mathcal{F}_0)\|_{\infty} \notag & \\
& \leq  2 \|\mathbb{E}(\|S_{2^{q-1}}\|_{\mathbb{H}}^2 \mid \mathcal{F}_0)\|_{\infty} +2 \| \mathbb{E}(<S_{2^{q-1}}-\mathbb{E}(S_{2^{q-1}} \mid \mathcal{F}_{2^{q-1}}),S_{2^q}-S_{2^{q-1}}>_{\mathbb{H}} \mid \mathcal{F}_0)\|_{\infty} \notag & \\
& \qquad \qquad + 2 \|\mathbb{E}(<\mathbb{E}(S_{2^{q-1}} \mid \mathcal{F}_{2^{q-1}}),S_{2^q}-S_{2^{q-1}}>_{\mathbb{H}} \mid \mathcal{F}_0) \|_{\infty}. \label{lemtech3} &
\end{align}
The last term in (\ref{lemtech3}) can be treated as in the proof of the corresponding facts in Proposition 2.1 of Peligrad and Utev \cite{peligradutev}, if we replace
everywhere the product in $\mathbb{R}$ by $<.,.>_{\mathbb{H}}$, and the $\mathbf{L}^2$-norm $\|x\|$ by the infinite norm. Consequently, we derive
\begin{align}
& \|\mathbb{E}(<\mathbb{E}(S_{2^{q-1}} \mid \mathcal{F}_{2^{q-1}}),S_{2^q}-S_{2^{q-1}}>_{\mathbb{H}} \mid \mathcal{F}_0)\|_{\infty} \notag & \\
&  \leq \sqrt{\|\mathbb{E}(\|S_{2^{q-1}}\|_{\mathbb{H}}^2 \mid \mathcal{F}_0)\|_{\infty}} \  2^{(q-1)/2}(\Delta_q-\Delta_{q-1}). \label{lemtech5} &
\end{align}
In the same way, since $\|\|S_{2^{q-1}}-\mathbb{E}(S_{2^{q-1}}\mid \mathcal{F}_{2^{q-1}})\|_{\mathbb{H}}\|_{\infty}=2^{(q-1)/2}(\Delta'_q-\Delta'_{q-1})$, we have
\begin{align}
& \| \mathbb{E}(<S_{2^{q-1}}-\mathbb{E}(S_{2^{q-1}} \mid \mathcal{F}_{2^{q-1}}),S_{2^q}-S_{2^{q-1}}>_{\mathbb{H}} \mid \mathcal{F}_0) \|_{\infty} \notag & \\
& \leq \sqrt{\|\mathbb{E}(\|S_{2^{q-1}}-\mathbb{E}(S_{2^{q-1}} \mid \mathcal{F}_{2^{q-1}})\|_{\mathbb{H}}^2 \mid \mathcal{F}_0)\|_{\infty}} \sqrt{\|\mathbb{E}(\|S_{2^q}-S_{2^{q-1}}\|_{\mathbb{H}}^2 \mid \mathcal{F}_0)\|_{\infty}} \notag & \\
& \leq \|\|S_{2^{q-1}}-\mathbb{E}(S_{2^{q-1}} \mid \mathcal{F}_{2^{q-1}})\|_{\mathbb{H}} \|_{\infty} \sqrt{\|\mathbb{E}(\|S_{2^{q-1}} \|_{\mathbb{H}}^2 \mid \mathcal{F}_0)\|_{\infty}} \notag & \\
& \leq \sqrt{\|\mathbb{E}(\|S_{2^{q-1}}\|_{\mathbb{H}}^2 \mid \mathcal{F}_0)\|_{\infty}} \ 2^{(q-1)/2}(\Delta'_q - \Delta'_{q-1}). \label{lemtech6} &
\end{align}
By induction and combining (\ref{lemtech5}) and (\ref{lemtech6}), we conclude that
\begin{align*}
& \quad \| \mathbb{E}(\|S_{2^q}\|_{\mathbb{H}} ^2 \mid \mathcal{F}_0) \|_{\infty}  \notag &  \\
& \leq  2 \times 2^{q-1}( \|\mathbb{E}(\|X_1\|_{\mathbb{H}}^2 \mid \mathcal{F}_0)\|_{\infty} + \frac{1}{2} \Delta_{q-1}+\frac{1}{2} \Delta'_{q-1})^2  \notag &  \\
& \qquad + 2 \times 2^{(q-1)/2}(\|\mathbb{E}(\|X_1\|_{\mathbb{H}}^2 \mid \mathcal{F}_0)\|_{\infty} + \frac{1}{2} \Delta_{q-1} + \frac{1}{2} \Delta'_{q-1}) \times 2^{(q-1)/2}(\Delta_q-\Delta_{q-1})  \notag &  \\
& \qquad +2 \times 2^{(q-1)/2}( \|\mathbb{E}(\|X_1 \|_{\mathbb{H}}^2 \mid \mathcal{F}_0)\|_{\infty} + \frac{1}{2} \Delta_{q-1}+\frac{1}{2} \Delta'_{q-1}) \times 2^{(q-1)/2}(\Delta'_q-\Delta'_{q-1})  \notag  & \\
& \leq  2^q(\|\mathbb{E}(\|X_1 \|_{\mathbb{H}}^2 \mid \mathcal{F}_0)\|_{\infty} +\frac{1}{2}\Delta_q + \frac{1}{2} \Delta'_q)^2  \notag & \\
& \leq   n (\|\mathbb{E}(\|X_1 \|_{\mathbb{H}}^2 \mid \mathcal{F}_0)\|_{\infty} +\frac{1}{2}\Delta_q + \frac{1}{2} \Delta'_q)^2.  \label{lemtech7} &
\end{align*}
\hfill $\square$
\subsubsection{Proof of Corollary \ref{corollaire}}
$\newline$
\indent The proof splits in two parts, and uses the same arguments as in the proof of Lemma $28$ in Dedecker, Merlev\`ede, Peligrad and Utev \cite{dedeckermerlevedepeligradutev}.
\begin{lem} \label{lemcor}
Assume that $\|\|X_0\|_{\mathbb{H}} \|_{\infty} < \infty $.
\begin{enumerate}
\item [i.] Under (\ref{conditionfonctionnel}) and (\ref{eqnfonct1cor}), (\ref{eqnfonct1}) holds. \\
\item [ii.] Under (\ref{conditionfonctionnel}) and (\ref{eqnfonct2cor}), (\ref{eqnfonct2})  holds.
\end{enumerate}
\end{lem}
\noindent \textit{Proof of Lemma \ref{lemcor}}.The proofs of $(i)$ and $(ii)$ are quite similarly, so here we prove only $(ii)$.\\
\indent Firstly, as in the proof of Lemma $28$ in Dedecker, Merlev\`ede, Peligrad and Utev \cite{dedeckermerlevedepeligradutev}, we prove by diadic recurrence (\ref{eqnfonct2}). Let $S_{a,b}=S_b-S_a$. Denote,
for any $t $ integer,
\begin{equation*}
A_{t,k}=\| \mathbb{E}(\|S_t\|_{\mathbb{H}}^2 \mid \mathcal{F}_{-k})- \mathbb{E}(\|S_t\|_{\mathbb{H}}^2)\|_{\infty}.
\end{equation*}
By stationarity, we have
\begin{eqnarray*}
A_{2t,k} & = & \| \mathbb{E}(\|S_{2t}\|^2_{\mathbb{H}} \mid \mathcal{F}_{-k})-\mathbb{E}(\|S_{2t}\|_{\mathbb{H}}^2)\|_{\infty} \notag \\
         & \leq & 2 \| \mathbb{E}(\|S_t \|_{\mathbb{H}}^2 \mid \mathcal{F}_{-k})- \mathbb{E}(\|S_t\|_{\mathbb{H}}^2)\|_{\infty} \notag \\
         &      & \quad +2 \| \mathbb{E}(<S_t,S_{t,2t}>_{\mathbb{H}} \mid \mathcal{F}_{-k})\|_{\infty} +2 \mid \mathbb{E}(<S_t,S_{t,2t}>_{\mathbb{H}}) \mid. \label{eqncor9}
\end{eqnarray*}
Moreover by Cauchy-Schwarz inequality and Lemma $\ref{lemtech}$, we get that
\begin{eqnarray*}
A_{2t,k} & \leq & 2 A_{t,k} + 4 \sqrt{\|\mathbb{E}(\|S_t\|_{\mathbb{H}}^2 \mid \mathcal{F}_0)\|_{\infty}} \|\|\mathbb{E}(S_t \mid \mathcal{F}_0)\|_{\mathbb{H}} \|_{\infty} \notag \\
 & & \qquad  + 4 \|\|S_t-\mathbb{E}(S_t \mid \mathcal{F}_t)\|_{\mathbb{H}}\|_{\infty} \sqrt{\|\mathbb{E}(\|S_t\|_{\mathbb{H}}^2 \mid \mathcal{F}_0)\|_{\infty}} \notag  \\
         & \leq & 2 A_{t,k} + 4 t^{1/2} \Delta_{\infty} \big \{ \|\|\mathbb{E}(S_t \mid \mathcal{F}_0)\|_{\mathbb{H}}\|_{\infty} + \|\|S_t-\mathbb{E}(S_t \mid \mathcal{F}_t)\|_{\mathbb{H}}\|_{\infty} \big \}. \notag \\
         & &   \label{eqncor10}
\end{eqnarray*}
With the notation
\begin{equation*}
B_{r,k}=2^{-r} \|\mathbb{E}(\|S_{2^r} \|_{\mathbb{H}}^2 \mid \mathcal{F}_{-k})-\mathbb{E}(\|S_{2^r}\|_{\mathbb{H}}^2) \|_{\infty}=2^{-r}A_{2^r,k},
\end{equation*}
by recurrence, for all $r \geq m$ and all $k >0$, we derive
\begin{eqnarray*}
B_{r,k} & \leq &  B_{r-1,k} + 2^{\frac{-r+3}{2} } \Delta_{\infty} \{ \|\|S_{2^{r-1}}-\mathbb{E}(S_{2^{r-1}} \mid \mathcal{F}_{2^{r-1}})\|_{\mathbb{H}}\|_{\infty} + \|\|\mathbb{E}(S_{2^{r-1}} \mid \mathcal{F}_0)\|_{\mathbb{H}} \|_{\infty} \} \\
& \leq & B_{m,k} +2 \Delta_{\infty} \Big \{ \sum_{j=m}^r 2^{-j/2} \|\|S_{2^j} - \mathbb{E}(S_{2^j} \mid \mathcal{F}_{2^j})\|_{\mathbb{H}}\|_{\infty} + \sum_{j=m}^r 2^{-j/2} \|\|\mathbb{E}(S_{2^j} \mid \mathcal{F}_0)\|_{\mathbb{H}} \|_{\infty} \Big \} \\
 & \leq & B_{m,k} +2 \Delta_{\infty} \{ \Delta_{m,\infty} + \Delta'_{m,\infty} \},
\end{eqnarray*}
where
\begin{equation*}
\Delta_{m,\infty}=\sum_{j=m}^{\infty} 2^{-j/2} \|\|\mathbb{E}(S_{2^j} \mid \mathcal{F}_0)\|_{\mathbb{H}}\|_{\infty} \ \mbox{and} \ \Delta'_{m,\infty}=\sum_{j=m}^{\infty} 2^{-j/2}\|\|S_{2^j}-\mathbb{E}(S_{2^j} \mid \mathcal{F}_{2^j})\|_{\mathbb{H}}\|_{\infty}.
\end{equation*}
By stationarity and triangle inequality,
\begin{align}
& \| \mathbb{E}(\|S_{2^r}\|_{\mathbb{H}}^2 \mid \mathcal{F}_0)-\mathbb{E}(\|S_{2^r}\|_{\mathbb{H}}^2)\|_{\infty} \notag & \\
 \leq & \| \mathbb{E}(\|S_{2^r}\|^2_{\mathbb{H}} - \|S_{k,k+2^r} \|_{\mathbb{H}}^2 \mid \mathcal{F}_0)\|_{\infty} + \| \mathbb{E}(\|S_{2^r}\|_{\mathbb{H}}^2 \mid \mathcal{F}_{-k})-\mathbb{E}(\|S_{2^r}\|_{\mathbb{H}}^2)\|_{\infty}, \notag &
\end{align}
we have, for all integer $r \geq m+1$,
\begin{align}
& 2^{-r} \| \mathbb{E}(\|S_{2^r} \|_{\mathbb{H}}^2 \mid \mathcal{F}_0) - \mathbb{E}(\|S_{2^r}\|_{\mathbb{H}}^2)\|_{\infty} \notag & \\
\leq & B_{m,k} +2 \Delta_{\infty} ( \Delta_{m,\infty} + \Delta'_{m,\infty}) +2^{-r/2+2} k \| \mathbb{E}(\|X_1\|_{\mathbb{H}}^2 \mid \mathcal{F}_0)\|_{\infty}^{1/2} \Delta_{\infty}. \notag &
\end{align}
Consequently,
\begin{align*}
& \underset{ r \longrightarrow \infty}{\limsup} \  2^{-r} \| \mathbb{E}(\|S_{2^r}\|^2_{\mathbb{H}} \mid \mathcal{F}_0)- \mathbb{E}(\|S_{2^r}\|^2_{\mathbb{H}})\|_{\infty} \notag  & \\
 & \qquad  \leq   B_{m,k} +2 \Delta_{\infty} (\Delta_{m,\infty} + \Delta'_{m,\infty}). \label{eqncor11} &
\end{align*}
Letting $k \rightarrow \infty$, and using condition (\ref{eqnfonct2cor}), it follows that $\underset{ k \longrightarrow \infty}{\lim} B_{m,k} =0$. Next letting $ m \rightarrow \infty$, and using condition (\ref{conditionfonctionnel}),
we then derive that
\begin{equation} \label{eqncor12}
\underset{ r \longrightarrow \infty}{\lim} \ 2^{-r} \| \mathbb{E}(\|S_{2^r}\|^2_{\mathbb{H}} \mid \mathcal{F}_0)-\mathbb{E}(\|S_{2^r}\|^2_{\mathbb{H}}) \|_{\infty}=0.
\end{equation}
To finish the proof, we use the diadic expansion $n=\sum_{k=0}^{r-1} 2^k a_k$, where $a_{r-1}=1$ and $a_k \in \{ 0,1 \} $, as the proof of Proposition $2.1$ in Peligrad and Utev \cite{peligradutev} in order to treat
the whole sequence $S_n$, for $2^{r-1} \leq n < 2^r$. We then use the following representation,
\begin{equation*}
S_n=\sum_{j=0}^{r-1} T_{2^j} a_j \ \mbox{where} \ T_{2^j}=\sum_{i=n_{j-1}+1}^{n_j} X_i, \ n_j=\sum_{k=0}^j 2^k a_k, \ n_{-1}=0.
\end{equation*}
Notice that
\begin{align}
& \frac{1}{n} \| \mathbb{E}(\|S_n\|^2_{\mathbb{H}} \mid \mathcal{F}_0)-\mathbb{E}(\|S_n\|_{\mathbb{H}}^2)\|_{\infty} \notag \\
 \leq & \frac{1}{n} \sum_{j=0}^{r-1} a_j \|\mathbb{E}( \|S_{2^j}\|_{\mathbb{H}}^2 \mid \mathcal{F}_0)-\mathbb{E}(\|S_{2^j}\|_{\mathbb{H}}^2)\|_{\infty} \notag \\
       &  +\frac{1}{n} \sum_{i \not=j= 0}^{r-1} a_i a_j \| \mathbb{E}(<T_{2^i},T_{2^j}>_{\mathbb{H}} \mid \mathcal{F}_0)-\mathbb{E}(<T_{2^i},T_{2^j}>_{\mathbb{H}}) \|_{\infty}. \label{eqnnonadapt1}
\end{align}
For the first term of the right-hand side in $(\ref{eqnnonadapt1})$, we treat it as a diadic integer,
\begin{equation*}\label{eqnnonadapt2}
\underset{n \longrightarrow \infty}{\lim} \frac{1}{n} \sum_{j=0}^{r-1} a_j \| \mathbb{E}(\|S_{2^j}\|_{\mathbb{H}}^2 \mid \mathcal{F}_0)-\mathbb{E}(\|S_{2^j}\|_{\mathbb{H}}^2)\|_{\infty}=0.
\end{equation*}
Suppose that $i<j<r$, we then have
\begin{align}
& \|\mathbb{E}(<T_{2^j},T_{2^i}>_{\mathbb{H}} \mid \mathcal{F}_0)-\mathbb{E}(<T_{2^j},T_{2^i}>_{\mathbb{H}})\|_{\infty} \notag \\
  \leq & \|\mathbb{E}(<T_{2^i}-\mathbb{E}(T_{2^i} \mid \mathcal{F}_{n_i}),T_{2^j}>_{\mathbb{H}} \mid \mathcal{F}_0)-\mathbb{E}(<T_{2^i}-\mathbb{E}(T_{2^i} \mid \mathcal{F}_{n_i}),T_{2^j}>_{\mathbb{H}})\|_{\infty} \notag \\
       &  + \| \mathbb{E}(<\mathbb{E}(T_{2^i} \mid \mathcal{F}_{n_i}),T_{2^j}>_{\mathbb{H}} \mid \mathcal{F}_0)-\mathbb{E}(<\mathbb{E}(T_{2^i} \mid \mathcal{F}_{n_i}),T_{2^j}>_{\mathbb{H}})\|_{\infty}. \label{eqnnonadapt3}
\end{align}
For the first term in the right-hand side in $(\ref{eqnnonadapt3})$, we have by Cauchy-Schwarz inequality,
\begin{align*}
& \sum_{i=0}^{r-2} \sum_{j=i+1}^{r-1} \| \mathbb{E}(<T_{2^i}-\mathbb{E}(T_{2^i} \mid \mathcal{F}_{n_i}),T_{2^j}>_{\mathbb{H}} \mid \mathcal{F}_0) \|_{\infty}\\
\leq &  \sum_{i=0}^{r-2} \sum_{j=i+1}^{r-1} \|\|T_{2^i}-\mathbb{E}(T_{2^i} \mid \mathcal{F}_{n_i})\|_{\mathbb{H}}\|_{\infty} \|\mathbb{E}(\|T_{2^j}\|_{\mathbb{H}}^2 \mid \mathcal{F}_0)\|_{\infty}^{1/2}.
\end{align*}
By $(\ref{eqncor12})$,
\begin{equation*}
\| \mathbb{E}(\|S_{2^r}\|_{\mathbb{H}}^2 \mid \mathcal{F}_0)\|_{\infty}^{1/2} = O(2^{r/2}).
\end{equation*}
Hence, we get
\begin{align*}
& \sum_{i=0}^{r-2} \sum_{j=i+1}^{r-1} \| \mathbb{E}(<T_{2^i}-\mathbb{E}(T_{2^i} \mid \mathcal{F}_{n_i}),T_{2^j}>_{\mathbb{H}} \mid \mathcal{F}_0) \|_{\infty} \\
\leq & C_1 \sum_{i=0}^{r-2} \|\|T_{2^i}-\mathbb{E}(T_{2^i} \mid \mathcal{F}_{n_i})\|_{\mathbb{H}}\|_{\infty} \sum_{j=i+1}^{r-1} 2^{j/2} \\
\leq & C_2 \sum_{i=0}^{r-2} \|\|T_{2^i}-\mathbb{E}(T_{2^i} \mid \mathcal{F}_{n_i})\|_{\mathbb{H}}\|_{\infty} 2^{r/2},
\end{align*}
where $C_1$ and $C_2$ are positive constants. \\
\indent Therefore, for all $2^{r-1} \leq n < 2^r$, we obtain
\begin{align*}
& \frac{1}{n} \sum_{i=0}^{r-2} \sum_{j=i+1}^{r-1} \| \mathbb{E}(<T_{2^i}-\mathbb{E}(T_{2^i} \mid \mathcal{F}_{n_i}),T_{2^j}>_{\mathbb{H}} \mid \mathcal{F}_0)\|_{\infty} \\
\leq &  C_2 2^{-r/2+1} \sum_{i=0}^{r-2} 2^{i/2} \frac{ \|\|S_{2^i}-\mathbb{E}(S_{2^i} \mid \mathcal{F}_{2^i})\|_{\mathbb{H}}\|_{\infty}}{2^{i/2}}.
\end{align*}
As
\begin{equation*}
\sum_{i=0}^{\infty} \frac{\|\|S_{2^i}-\mathbb{E}(S_{2^i} \mid \mathcal{F}_{2^i})\|_{\mathbb{H}}\|_{\infty}}{2^{i/2}} < \infty,
\end{equation*}
we conclude by Kronecker lemma that
\begin{equation}\label{eqnnonadapt4}
\underset{n\longrightarrow \infty}{\lim} \frac{1}{n} \sum_{i=0}^{r-2} \sum_{j=i+1}^{r-1} \|\mathbb{E}(<T_{2^i}-\mathbb{E}(T_{2^i}\mid \mathcal{F}_{n_i}),T_{2^j}>_{\mathbb{H}} \mid \mathcal{F}_0)\|_{\infty}=0.
\end{equation}
For the second term of the right-hand side in $(\ref{eqnnonadapt3})$, we use the same arguments as in the proof of Proposition 2.1 in Peligrad and Utev \cite{peligradutev} by replacing
the product in $\mathbb{R}$, by $<.,.>_{\mathbb{H}}$ and the $\mathbf{L}^2$-norm by the infinite norm.
Consequently, we get
\begin{equation} \label{eqnnonadapt5}
\underset{n \longrightarrow \infty}{\lim} \frac{1}{n} \sum_{0 \leq i < j \leq r-1} \| \mathbb{E}(<\mathbb{E}(T_{2^i} \mid \mathcal{F}_{n_i}),T_{2^j}>_{\mathbb{H}} \mid \mathcal{F}_0)-\mathbb{E}(< \mathbb{E}(T_{2^i} \mid \mathcal{F}_{n_i}),T_{2^j}>_{\mathbb{H}})\|_{\infty}=0.
\end{equation}
Combining $(\ref{eqnnonadapt4})$ and $(\ref{eqnnonadapt5})$, we conclude
\begin{equation*}
\underset{n \longrightarrow \infty}{\lim} \frac{1}{n} \sum_{i \not= j=0}^{r-1} a_i a_j \| \mathbb{E}(<T_{2^i},T_{2^j}>_{\mathbb{H}} \mid \mathcal{F}_0)-\mathbb{E}(<T_{2^i},T_{2^j}>_{\mathbb{H}})\|_{\infty}=0.
\end{equation*}
This proves $(\ref{eqnfonct2})$.
\hfill $\square$
\subsection{Proof of  Proposition \ref{propfunclinearprocess}}
$\newline$
\indent Let $\varepsilon'$ be an independent copy of $\varepsilon$, and denote by $\mathbb{E}_{\varepsilon}(.)$ the conditional expectation
with respect to $\varepsilon$. Define
\begin{equation*}
Y_n=\sum_{i<n} c_i(\varepsilon_{n-i}), Y_n'=\sum_{i<n} c_i(\varepsilon'_{n-i}), Z_n=\sum_{i \geq n} c_i(\varepsilon_{n-i}),Z_n'=\sum_{i \geq n} c_i(\varepsilon'_{n-i}).
\end{equation*}
Then, taking $\mathcal{F}_l=\sigma(\varepsilon_i,i \leq l)$, we have
\begin{eqnarray*}
\|\|\mathbb{E}(X_n \mid \mathcal{F}_0)\|_{\mathbb{H}}\|_{\infty} &= & \|\|\mathbb{E}_{\varepsilon}[f(Y'_n+Z_n)-f(Y'_n+Z'_n)] \|_{\mathbb{H}}\|_{\infty} \\
                                                                 & \leq & w_f \big(\|\|\varepsilon_0-\varepsilon'_0\|_{\mathbb{H}}\|_{\infty} \sum_{k \geq n} \|c_k\|_{L(\mathbb{H})} \big)
\end{eqnarray*}
and
\begin{eqnarray*}
\|\|X_{-n}-\mathbb{E}(X_{-n} \mid \mathcal{F}_0)\|_{\mathbb{H}}\|_{\infty} & = & \|\|\mathbb{E}_{\varepsilon}(f(Y_{-n}+Z_{-n})-f(Z_{-n}+Y'_{-n}))\|_{\mathbb{H}}\|_{\infty} \\
                                                                           & \leq & w_{f} \big( \|\|\varepsilon_0-\varepsilon'_0\|_{\mathbb{H}}\|_{\infty} \sum_{k \leq -n} \|c_k\|_{L(\mathbb{H})} \big).
\end{eqnarray*}
Then the condition (\ref{conditionfonctionnelcor}) is satisfied as soon as (\ref{eqn427}) holds. \\
\indent As the proof of (\ref{eqnfonct1cor}) is quite similar of the proof of (\ref{eqnfonct2cor}), we only prove (\ref{eqnfonct2cor}).\\
\indent We have for all  integer $p \geq 0$, \\
$\|  \mathbb{E}(<X_i,X_{i+p}>_{\mathbb{H}} \mid \mathcal{F}_0)-\mathbb{E}(<X_i,X_{i+p}>_{\mathbb{H}})\|_{\infty} \notag $
\begin{eqnarray}
 &\leq & C \Big \{ w_f \Big( \| \| \varepsilon_0-\varepsilon'_0 \|_{\mathbb{H}} \|_{\infty} \sum_{k \geq i} \|c_k\|_{L(\mathbb{H})} \Big) \notag \\
  && \qquad \qquad +w_f \Big( \| \| \varepsilon_0-\varepsilon'_0 \|_{\mathbb{H}} \|_{\infty} \sum_{k \geq i+p} \|c_k\|_{L (\mathbb{H})} \Big) \Big \}, \label{eqn30007}
\end{eqnarray}
where $C$ is a constant.\\
\indent By $(\ref{eqn427})$ and Corollary \ref{corollaire}, Proposition \ref{propfunclinearprocess} holds.
\hfill $\square$
\subsection{Proof of Proposition \ref{propstablemarkovchain}}
$\newline$
\indent Firstly, we give a technical lemma,
\begin{lem} \label{lemstablemarkovchain}
If $ \mathrm{Lip}(K^n(f)) \leq C \rho^n \mathrm{Lip}(f)$, then
\begin{equation} \label{eqnchainemarkov1}
\|\mathbb{E}(f(Y_k)\mid Y_0)-\mathbb{E}(f(Y_k))\|_{\infty} \leq 2 \|\|Y_0\|_{\mathbb{H}}\|_{\infty} C \rho^k \mathrm{Lip}(f).
\end{equation}
\end{lem}
\noindent \textit{Proof of Lemma \ref{lemstablemarkovchain}}.
As
\begin{equation*}
\mathbb{E}(f(Y_k)\mid Y_0=y)-\mathbb{E}(f(Y_k)) = \int \big (K^k(f)(y)-K^k(f)(z) \big) \mu(dz),
\end{equation*}
we deduce
\begin{eqnarray}
\|\|\mathbb{E}(f(Y_k)\mid Y_0=y)-\mathbb{E}(f(Y_k))\|_{\mathbb{H}}\|_{\infty} & \leq & \Big \|\int \|K^k(f)(y)-K^k(f)(z)\|_{\mathbb{H}} \mu(dz) \Big \|_{\infty} \notag \\
                                                                              & \leq & \mbox{Lip}(K^k(f)) \ \Big \|\int \|y-z\|_{\mathbb{H}} \mu(dz) \Big \|_{\infty} \notag \\
                                                                              & \leq & C \rho^k \mbox{Lip}(f) \ \Big \|\int \|y-z\|_{\mathbb{H}} \mu (dz) \Big \|_{\infty}. \notag \\
                                                                              & &  \label{eqnchainemarkov2}
\end{eqnarray}
Observe that
\begin{equation} \label{eqnchainemarkov3}
\Big \|\int \|y-z\|_{\mathbb{H}} \mu(dz) \Big \|_{\infty} \leq \Big \|\int \big (\|y\|_{\mathbb{H}} + \|z \|_{\mathbb{H}} \big )   \mu(dz) \Big \|_{\infty} \leq 2 \|\|Y_0\|_{\mathbb{H}}\|_{\infty}.
\end{equation}
Consequently, combining (\ref{eqnchainemarkov2}) and (\ref{eqnchainemarkov3}), we have
\begin{equation*}
\|\|\mathbb{E}(f(Y_k)\mid Y_0=y)-\mathbb{E}(f(Y_k))\|_{\mathbb{H}}\|_{\infty} \leq 2 \|\|Y_0\|_{\mathbb{H}}\|_{\infty} C \rho^k \mbox{Lip}(f).
\end{equation*}
\hfill $\square$ \\
\noindent \textit{Proof of Proposition \ref{propstablemarkovchain}}. We apply Corollary \ref{corollaire} to the following random variables,
\begin{equation*}
X_k=f(Y_k)-\mathbb{E}(f(Y_k)), \ \forall \ k \geq 0.
\end{equation*}
Since $(Y_n)_{n \geq 0}$ is a Markov chain, we have to prove that
\begin{equation*} \label{eqnchainemarkov4}
\sum_{k \geq 1} \frac{1}{\sqrt{k}} \|\|\mathbb{E}(X_k \mid \mathcal{F}_0)\|_{\mathbb{H}}\|_{\infty}< \infty.
\end{equation*}
By Lemma \ref{lemstablemarkovchain}, we derive
\begin{equation*}
\sum_{k \geq 1} \frac{1}{\sqrt{k}} \|\|\mathbb{E}(X_k \mid \mathcal{F}_0)\|_{\mathbb{H}}\|_{\infty} \leq 2 \|\|Y_0\|_{\mathbb{H}}\|_{\infty} C \mbox{Lip}(f) \sum_{k \geq 1} \frac{1}{\sqrt{k}} \rho^k < \infty.
\end{equation*}
The proof of $(\ref{eqnfonct1cor})$ is quite similar of the proof of $(\ref{eqnfonct2cor})$, so we only detail $(\ref{eqnfonct2cor})$.
If $k>l$, by triangle inequality
\begin{align*}
& \|\mathbb{E}(<X_k,X_l>_{\mathbb{H}}\mid \mathcal{F}_{-n})-\mathbb{E}(<X_k,X_l>_{\mathbb{H}})\|_{\infty} \notag \\
\leq & \|\mathbb{E}(<f(Y_k),f(Y_l)>_{\mathbb{H}} \mid \mathcal{F}_{-n})-\mathbb{E}(<f(Y_k),f(Y_l)>_{\mathbb{H}})\|_{\infty} \notag \\
&  +\| \mathbb{E}(<f(Y_k), \mathbb{E}(f(Y_l))>_{\mathbb{H}} \mid \mathcal{F}_{-n})-\mathbb{E}(<f(Y_k),\mathbb{E}(f(Y_l))>_{\mathbb{H}})\|_{\infty} \notag \\
&  + \|\mathbb{E}(<\mathbb{E}(f(Y_k)),f(Y_l)>_{\mathbb{H}} \mid \mathcal{F}_{-n})-\mathbb{E}(<\mathbb{E}(f(Y_k)),f(Y_l)>_{\mathbb{H}})\|_{\infty}.
\end{align*}
Using Lemma \ref{lemstablemarkovchain}, we get
\begin{align*}
& \|\mathbb{E}(<f(Y_k),\mathbb{E}(f(Y_l))>_{\mathbb{H}} \mid \mathcal{F}_{-n})-\mathbb{E}(<f(Y_k),\mathbb{E}(f(Y_l))>_{\mathbb{H}})\|_{\infty} \\
\leq & \|\|\mathbb{E}(f(Y_l))\|_{\mathbb{H}}\|_{\infty} \|\|\mathbb{E}(f(Y_k) \mid \mathcal{F}_{-n})-\mathbb{E}(f(Y_k))\|_{\mathbb{H}}\|_{\infty} \\
\leq & 2 C \|\|\mathbb{E}(f(Y_0))\|_{\mathbb{H}}\|_{\infty} \|\|Y_0\|_{\mathbb{H}}\|_{\infty} \mbox{Lip}(f) \rho^{k+n} \underset{n \rightarrow \infty}{\longrightarrow} 0,
\end{align*}
and for $k>l$,
\nolinebreak
\begin{align*}
& \|\mathbb{E}(<f(Y_k),f(Y_l)>_{\mathbb{H}} \mid \mathcal{F}_{-n})-\mathbb{E}(<f(Y_k),f(Y_l)>_{\mathbb{H}})\|_{\infty} \\
= & \|\mathbb{E}(<\mathbb{E}(f(Y_k) \mid \mathcal{F}_{l}),f(Y_l)>_{\mathbb{H}} \mid \mathcal{F}_{-n})-\mathbb{E}(<\mathbb{E}(f(Y_k) \mid \mathcal{F}_{l}),f(Y_l)>_{\mathbb{H}})\|_{\infty} \\
= & \|\mathbb{E}(<K^{k-l}(f)(Y_l),f(Y_l)>_{\mathbb{H}} \mid \mathcal{F}_{-n})-\mathbb{E}(<K^{k-l}(f)(Y_l),f(Y_l)>_{\mathbb{H}})\|_{\infty} \\
= & \| K^{l+n}(<K^{k-l}(f)(.),f(.)>_{\mathbb{H}})(Y_{-n})-\mu(K^{l+n}(<K^{k-l}(f)(.),f(.)>_{\mathbb{H}})(Y_{-n}))\|_{\infty} \\
\leq & 2 C \|\|Y_0\|_{\mathbb{H}}\|_{\infty} \mbox{Lip}(<K^{k-l}(f)(.),f(.)>_{\mathbb{H}}) \rho^{l+n} \ \underset{n \rightarrow \infty}{\longrightarrow} 0.
\end{align*}
\nolinebreak
Hence (\ref{eqnfonct2cor}) holds.
\hfill $\square$
\subsection{Proof of Proposition \ref{propcramervonmises}}
$\newline$
\indent We apply Corollary \ref{corollaire} to the random variables $X_i=\{t \mapsto \mathbf{1}_{Y_i\leq t}-\mathbb{F}(t): \ t \in \mathbb{R} \}$. Since
\begin{equation*}
\sum_{n \geq 1} \frac{1}{\sqrt{n}} \tilde{\phi_2}(n) < \infty \Longrightarrow \sum_{n \geq 1} \frac{1}{\sqrt{n}} \tilde{\phi_1}(n) < \infty,
\end{equation*}
the condition (\ref{conditionfonctionnelcor}) holds. \\
\indent As the proofs of (\ref{eqnfonct1cor}) and (\ref{eqnfonct2cor}) are quite similar, we only detail the proof of (\ref{eqnfonct2cor}).\\
\indent By Fubini, we have, for any $i<j$,
\begin{align*}
&\|\mathbb{E}(<X_i,X_j>_{\mathbf{L}^2(\mathbb{R},\mu)} \mid \mathcal{F}_0)-\mathbb{E}(<X_i,X_j>_{\mathbf{L}^2(\mathbb{R},\mu)})\|_{\infty} \notag & \\
= & \Big \| \mathbb{E} \Big( \int(\mathbf{1}_{Y_i \leq t }-\mathbb{F}(t))(\mathbf{1}_{Y_j \leq t} - \mathbb{F}(t)) \, \mu(dt) \Big | \mathcal{F}_0 \Big) \notag \\
& \qquad \qquad -\mathbb{E} \Big (\int ( \mathbf{1}_{Y_i \leq t}-\mathbb{F}(t))(\mathbf{1}_{Y_j \leq t}-\mathbb{F}(t)) \, \mu(dt) \Big ) \Big \|_{\infty} \notag & \\
\leq & \int \| \mathbb{E}((\mathbf{1}_{Y_i \leq t } -\mathbb{F}(t))(\mathbf{1}_{Y_j \leq t}-\mathbb{F}(t)) \mid \mathcal{F}_0)-\mathbb{E}((\mathbf{1}_{Y_i \leq t}-\mathbb{F}(t))( \mathbf{1}_{Y_j \leq t}-\mathbb{F}(t)))\|_{\infty} \, \mu(dt) \notag & \\
\leq & \|b(\mathcal{F}_0,Y_i,Y_j)\|_{\infty} \leq \tilde{\phi_2}(i). \label{eqncvm2} &
\end{align*}
Since $\sum_{n \geq 1} n^{-1/2} \tilde{\phi_2}(i)<\infty, \ \tilde{\phi_2}(i) \underset{i \longrightarrow \infty}{\longrightarrow} 0$, all
the conditions of Corollary \ref{corollaire} are true. \\
\indent From Dedecker and Merlev\`ede \cite{dedeckermerlevedeTCL}, the $\mathbf{L}^2(\mathbb{R},\mu)$-valued random variable $\sqrt{n} (\mathbb{F}_n-\mathbb{F})$ converges
stably to a zero mean $\mathbf{L}^2(\mathbb{R},\mu)$-valued gaussian random variable $\mathbb{G}$, with covariance function $Q$, given in Proposition \ref{propcramervonmises}. \\
\indent We deduce that $\sqrt{n} (\mathbb{F}_n-\mathbb{F})$ satisfies the MDP in $\mathbf{L}^2(\mathbb{R},\mu)$, with the good rate function
\begin{equation*} \label{eqncvm3}
\forall \ f \in \mathbf{L}^2(\mathbb{R},\mu), \ I(f)=\underset{g \in \mathbf{L}^2(\mathbb{R},\mu)} {\sup} (<f,g>_{\mathbf{L}^2(\mathbb{R},\mu)}-\frac{1}{2}<g,Qg>_{\mathbf{L}^2(\mathbb{R},\mu)}).
\end{equation*}
\hfill $\square$
$\nobreak$
\section{appendix}
\begin{lem}\label{appendix} Let $(U_j)_{j \geq 0}$ be a sequence of positive reals such that $U_0=0$ and $U_{i+j} \leq \tilde{C_1} U_i+ \tilde{C_2} U_j$. Let $C=\tilde{C_1}+\tilde{C_2}$. Then,
\begin{enumerate}
\item [$1.$] For $n$, and $r$ integers such that $n \geq 1$, $2^{r-1} \leq n < 2^r$, and $p \geq 1$,
\nolinebreak
\begin{equation} \label{eqnappendix1}
\sum_{j=0}^{r-1} \frac{1}{2^{j(p-1)}} U_{2^j} \leq \frac{C}{(1-2^{-p})} \sum_{k=1}^{n-1} \frac{1}{k^p} U_k.
\end{equation}
\nolinebreak
\item [$2.$] If $\sum_{k=1}^{\infty} k^{-p} U_k < \infty$ for a $p>1$, then
\nolinebreak
\begin{equation} \label{eqnappendix2}
\frac{1}{m^{p-1}} \sum_{j \geq 1} \frac{1}{j^p} U_{jm} \underset{ m \longrightarrow \infty}{\longrightarrow} 0.
\end{equation}
In particular,
\begin{equation*}
\frac{U_m}{m^{p-1}} \underset{ m \longrightarrow \infty}{\longrightarrow} 0.
\end{equation*}
\end{enumerate}
\end{lem}
\noindent \textit{Proof of Lemma \ref{appendix}}.
Firstly, we prove (\ref{eqnappendix1}). With this aim, we note that
\begin{equation*}
U_T \leq \tilde{C_1} U_k + \tilde{C_2} U_{T-k},
\end{equation*}
so
\begin{equation*}
(T+1) U_T \leq C \sum_{k=0}^T U_k.
\end{equation*}
Thus, for $n \geq 2^{r-1}$, we get
\begin{equation*}
\sum_{j=0}^{r-1} \frac{U_{2^j}}{2^{j(p-1)}} \leq C \sum_{k=1}^{2^{r-1}} U_k \sum_{j:2^j \geq k} \frac{1}{(2^j+1)2^{j(p-1)}} \leq \frac{C}{(1-2^{-p})} \sum_{k=1}^n \frac{U_k}{k^p}.
\end{equation*}
To prove (\ref{eqnappendix2}), we write
\begin{equation*}
\frac{1}{m^{p-1}} \sum_{j=1}^{\infty} \frac{U_{jm}}{j^p} \leq C \sum_{k=1}^{\infty} U_k \sum_{j:jm\geq k} \frac{1}{j^p(jm+1)m^{p-1}} \leq C c_p \sum_{k=1}^{\infty} \frac{U_k}{(k+m)^p} \underset{m\rightarrow \infty}{\longrightarrow} 0,
\end{equation*}
by the Fatou lemma since $U_k(k+m)^{-p}\downarrow 0$ as $m\rightarrow\infty$. Here, $c_p$ is a positive constant depending on $p$.
\hfill $\square$
\nobreak
\renewcommand{\refname}{References}

\end{document}